\documentclass[sn-mathphys,Numbered]{sn-jnl}
\usepackage[utf8]{inputenc}
\usepackage[table]{xcolor}

\newcommand{\domainBoundaryDirichlet}{\domainBoundary_{\deformationMap}}
\newcommand{\domainBoundaryTraction}{\domainBoundary_{\PKStress}}
\newcommand{\tractionBC}{\boldsymbol T_{\domainBoundary}}
\newcommand{\deformationMapBC}{\deformationMap_{\domainBoundary}}
\newcommand{\deformationGradient}{\boldsymbol F}
\newcommand{\nCoarseDofs}{\overline{N}_{\mathrm{dof}}}

\newcommand{\SPSDNN}{SPSD-NN}
\newcommand{\SPSDLLS}{SPSD-LLS}
\newcommand{\contactCaption}{Example 2, fastener contact.}
\newcommand{\preloadCaption}{Example 3, fastener with preload and contact.}
\newcommand{\cubeCaption}{Example 1, multi-element cube.}
\newcommand{\displacementScalar}{u}
\newcommand{\forcingScalar}{b}
\newcommand{\forceOffset}{\mathbf{f}_{\text{0}}}
\newcommand{\CauchyStress}{\boldsymbol \sigma}
\newcommand{\PKStress}{\boldsymbol P}

\newcommand{\PKStressFine}{\PKStress^{\prime}}

\newcommand{\referenceDomain}{\domain}

\newcommand{\deformedDomain}{\domain^*}

\newcommand{\materialCoordinate}{\boldsymbol X}
\newcommand{\eulerianCoordinate}{\boldsymbol x}
\newcommand{\deformationMap}{\boldsymbol \psi}
\newcommand{\deformationMapOne}{\overline{\deformationMap}}
\newcommand{\deformationMapTwo}{{\deformationMap}^{\prime}}

\newcommand{\deformationMapFEM}{\boldsymbol \psi_{\mathrm{h}}}
\newcommand{\deformationMapCoarseFEM}{\overline{\boldsymbol \psi}_{\mathrm{h}}}

\newcommand{\deformationMapCoarseDiscreteFEM}{\overline{\mathbf{\psi}}_\mathrm{h}}

\newcommand{\reducedMlModelVec}{\boldsymbol{\hat{\mlModel}}}

\newcommand{\reducedMlModelStiffness}{\hat{\mathcal{L}}}
\newcommand{\reducedMlModelStiffnessVec}{\boldsymbol{\hat{\mathcal{L}}}}

\newcommand{\closureTerm}{\overline{\mathbf{s}}}
\newcommand{\stiffnessMatrix}{\mathbf{K}}
\newcommand{\stiffnessMatrixCoarse}{\overline{\mathbf{K}}}
\newcommand{\reducedDimension}{K}
\newcommand{\reducedNeuralNetwork}{\boldsymbol{\mathcal{N}}}
\newcommand{\reducedNeuralNetworkStiffness}{\boldsymbol{\mathcal{L}}_{\mathcal{NN}}}

\newcommand{\weights}{\mathbf{w}}
\newcommand{\nWeights}{N_{\weights}}
\newcommand{\scaleFactorZOne}{b_{z1}}
\newcommand{\scaleFactorZTwo}{b_{z2}}
\newcommand{\timeTermZOne}{T_{z1}}
\newcommand{\cosRamp}[4]{\frac{1}{2}\left(#1 + #1 \cos\left( \left( \frac{t - #2}{#3 - #2}\right) \pi - \pi \right) \right)} 

\newcommand{\scaleFactorXOne}{b_{x1}}
\newcommand{\scaleFactorXTwo}{b_{x2}}
\newcommand{\timeTermXOne}{T_{x1}}

\newcommand{\domain}{\Omega}

\newcommand{\domainBoundary}{\Gamma}

\newcommand{\RR}[1]{\mathbb{R}^{#1}}
\newcommand{\nSamples}{N_s}
\newcommand{\forceSample}{\mathbf{F}}
\newcommand{\forceSampleOffset}{\mathbf{F}^*}

\newcommand{\domainBoundaryOneTwo}{\Gamma_{\mathcal{I}}}

\newcommand{\bz}{\mathbf{0}}

\newcommand{\normal}{\mathbf{n}}
\newcommand{\normalOneTwo}{\mathbf{n}_{\mathcal{I}}}

\newcommand{\domainOne}{\overline{\domain} }

\newcommand{\domainTwo}{\domain^{\prime}}
\newcommand{\forcing}{\boldsymbol b}

\newcommand{\mlModelVec}{\boldsymbol{ \mathcal{M}}}

\newcommand{\mlModel}{\mathcal{M}}

\newcommand{\reducedLinearRegressionMatrix}{\hat{\mathbf{A}}}

\newcommand{\basis}{\boldsymbol \Phi}
\newcommand{\basisCombined}{\boldsymbol \Phi_{*}}
\newcommand{\reducedDimensionCombined}{K_{*}}

\newcommand{\basisForces}{\boldsymbol \Phi_{\force}}
\newcommand{\basisDisplacements}{\boldsymbol \Phi_{\displacement}}
\newcommand{\reducedDimensionForces}{K_{\force}}
\newcommand{\reducedDimensionDisplacement}{K_{\displacement}}

\newcommand{\reducedLowerLinearRegressionMatrix}{\hat{\mathbf{L}}}

\newcommand{\force}{\boldsymbol{f}}

\newcommand{\displacement}{\boldsymbol u}
\newcommand{\nElements}{N_{\text{el}}}
\newcommand{\nCoarseElements}{\overline{N_{\text{el}}}}

\newcommand{\features}{\boldsymbol \theta}
\newcommand{\featuresMatrix}{\mathbf{U}}
\newcommand{\nFeatures}{N_{\features}}

\newcommand{\x}{\boldsymbol x}
\newcommand{\displacementFineDiscreteFEM}{{\mathbf{u}}^{\prime}_{\mathsf{h}}}
\newcommand{\displacementFineDiscreteScalarFEM}{{{u}}^{\prime}_{\mathsf{h}}}

\newcommand{\displacementCoarseDiscreteFEM}{\overline{\mathbf{u}}_{\mathsf{h}}}
\newcommand{\displacementCoarseDiscreteScalarFEM}{\overline{u}_{\mathsf{h}}}

\newcommand{\displacementCoarseDiscreteBoundaryFEM}{\overline{\mathbf{u}}^{\Gamma_\mathcal{I}}_{\mathsf{h}}}

\newcommand{\displacementDiscreteFEM}{\mathbf{u}_{\mathsf{h}}}
\newcommand{\displacementDiscreteScalarFEM}{{u}_{\mathsf{h}}}

\newcommand{\coarseNodeSet}{\overline{\mathcal{N}}}
\newcommand{\coarseNodeSetDomainBoundaryOneTwo}{\overline{\mathcal{N}}^{\domainBoundaryOneTwo}}

\newcommand{\nCoarseNodes}{\overline{N}}
\newcommand{\nCoarseNodesDomainBoundaryOneTwo}{\overline{N}^{\domainBoundaryOneTwo}}
\newcommand{\nCoarseDofsDomainBoundaryOneTwo}{\overline{N}^{\domainBoundaryOneTwo}_{\mathrm{dof}}}

\newcommand{\displacementFEM}{\boldsymbol{u}_{\mathrm{h}}}

\newcommand{\PKStressAsDeformationArg}[1]{\PKStress \left( #1  \right) }

\newcommand{\trialFunction}{\boldsymbol \phi}
\newcommand{\trialFunctionOne}{\overline{\trialFunction}}
\newcommand{\trialFunctionOneFEM}{{\trialFunctionOne_{\mathrm{h}}}}

\newcommand{\testFunction}{\boldsymbol v}
\newcommand{\testFunctionFEM}{\boldsymbol v_{\mathrm{h}}}

\newcommand{\testFunctionOne}{\overline{\testFunction}}

\newcommand{\testFunctionOneFEM}{\testFunctionOne_{\mathrm{h}}}

\newcommand{\coarseTrialSpaceFEM}{\overline{\trialSpaceFEM}}
\newcommand{\coarseTestSpaceFEM}{\overline{\testSpaceFEM}}

\newcommand{\trialSpace}{\mathcal{V}}

\newcommand{\trialSpaceOne}{\overline{\trialSpace}}

\newcommand{\testSpaceOne}{\overline{\testSpace}}

\newcommand{\testSpace}{\mathcal{W}}
\newcommand{\testSpaceFEM}{\mathcal{W}_{\mathrm{h}}}

\newcommand{\trialSpaceFEM}{\mathcal{V}_{\mathrm{h}}}

\newcommand{\svdLeft}{\boldsymbol U}
\newcommand{\svdLeftVecArg}[1]{\boldsymbol u_{#1}}
\newcommand{\svdMid}{\boldsymbol \Sigma}
\newcommand{\svdRight}{\boldsymbol V^T}

\usepackage{caption}[font=small,skip=-2pt]
\usepackage{subcaption}[font=small,skip=-2pt]
\usepackage{float}

\usepackage{amsthm}

\newtheorem{remark}{Remark}[section] 

\usepackage{comment}
\usepackage{graphicx}%
\usepackage{multirow}%
\usepackage{amsmath,amssymb,amsfonts}%
\usepackage{amsthm}%
\usepackage{mathrsfs}%
\usepackage[title]{appendix}%
\usepackage{xcolor}%
\usepackage{textcomp}%
\usepackage{manyfoot}%
\usepackage{booktabs}%
\usepackage{algorithm}%
\usepackage{algorithmicx}%
\usepackage{algpseudocode}%
\usepackage{listings}%

\usepackage[title]{appendix}%

\usepackage{nomencl}
  \makenomenclature

\title[Embedded symmetric positive semi-definite machine-learned elements for reduced-order modeling in finite-element simulations with application to threaded fasteners]{Embedded symmetric positive semi-definite machine-learned elements for reduced-order modeling in finite-element simulations with application to threaded fasteners}
\author*[1]{\fnm{Eric} \sur{Parish}}\email{ejparis@sandia.gov}
\author[2]{\fnm{Payton} \sur{Lindsay}}\email{plindsa@sandia.gov}
\author[2]{\fnm{Timothy} \sur{Shelton}}\email{trshelt@sandia.gov}
\author[2]{\fnm{John} \sur{Mersch}}\email{jpmersc@sandia.gov}

\affil*[1]{\orgname{Sandia National Laboratories}, \orgaddress{\street{7011 East Ave}, \city{Livermore}, \postcode{94550}, \state{CA}, \country{United States}}}

\affil*[2]{\orgname{Sandia National Laboratories}, \orgaddress{\street{1515 Eubank Blvd SE}, \city{Albuquerque}, \postcode{87123}, \state{NM}, \country{United States}}}

\abstract{We present a machine-learning strategy for finite element analysis of solid mechanics wherein we replace complex portions of a computational domain with a data-driven surrogate. In the proposed strategy, we decompose a computational domain into an ``outer" coarse-scale domain that we resolve using a finite element method (FEM) and an ``inner" fine-scale domain. We then develop a machine-learned (ML) model for the impact of the inner domain on the outer domain. In essence, for solid mechanics, our machine-learned surrogate performs static condensation of the inner domain degrees of freedom. This is achieved by learning the map from displacements on the inner-outer domain interface boundary to forces contributed by the inner domain to the outer domain on the same interface boundary. We consider two such mappings, one that directly maps from displacements to forces without constraints, and one that maps from displacements to forces by virtue of learning a symmetric positive semi-definite (SPSD) stiffness matrix. We demonstrate, in a simplified setting, that learning an SPSD stiffness matrix results in a coarse-scale problem that is well-posed with a unique solution. We present numerical experiments on several exemplars, ranging from finite deformations of a cube to finite deformations with contact of a fastener-bushing geometry. We demonstrate that enforcing an SPSD stiffness matrix drastically improves the robustness and accuracy of FEM--ML coupled simulations, and that the resulting methods can accurately characterize out-of-sample loading configurations with significant speedups over the standard FEM simulations.} 

\begin{document}

\maketitle

\section{Introduction}
\begin{wrapfigure}{r}{0.5\linewidth}
\begin{center}
\includegraphics[trim={0cm 0cm 0cm 0cm},clip,width=1.0\linewidth]{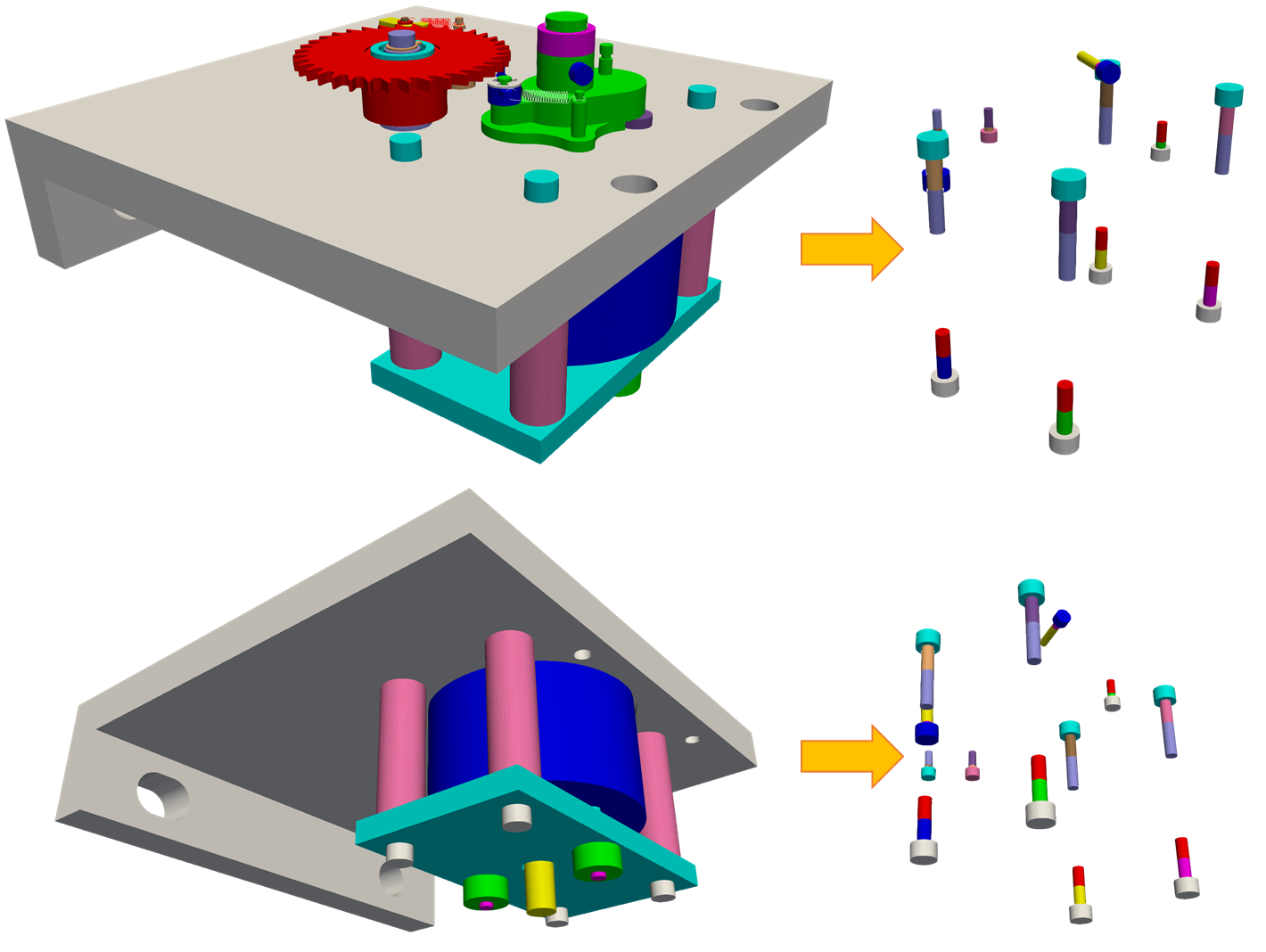}
\end{center}
\caption{Schematic of a ratcheting mechanism. Fasteners are an integral part of this complex mechanism and must capture relevant mechanics to properly assess quantities of interest.}
\label{fig:fastener_example}
\end{wrapfigure}

Finite element analysis of component and assembly-level systems remains computationally intensive despite the tremendous advancements in algorithmic research and computing power that have occurred over the past several decades. This problem is exacerbated when the underlying problem is multiscale in nature, wherein the physical discretization of the governing equations requires resolving a wide range of length and/or time scales. One such example that motivates the present work is that of systems-level models that involve threaded fasteners. In such systems, fasteners can be an integral connector of many sub-components, and accurately modeling the behavior of the fastener is critical for many analyses. Unfortunately, directly resolving the full fastener within a finite element model is difficult: the geometries are challenging to mesh, the resulting meshes can require hundreds of thousands of degrees of freedom, and the required mechanics that need to be simulated are complex as they often involve contact, friction, etc. As an example, Figure~\ref{fig:fastener_example} depicts a schematic of a ratcheting mechanism; fasteners are an essential part of this complex mechanism and must capture relevant mechanics to properly assess quantities of interest. In general, directly resolving each fastener is often not feasible in systems-level models, and approaches capable of reducing this computational burden are needed.


Various approaches have been developed to reduce computational complexities in finite element methods. These methods include the variational multiscale method~\cite{HuFeMa98}, Guyan reduction~\cite{Gu65}, Craig--Bampton reduction~\cite{CrBa68}, the generalized finite element method~\cite{FrBe18}, multiscale methods~\cite{HoWu97}, and reduced basis/proper-orthogonal decomposition reduced-order models~\cite{ChOlSe20,LiFiTe22,HoRoVe01,HoRoVe02}. Machine learning-based techniques have garnered significant attention recently and are the focus of the present literature review. ``Smart" finite elements were introduced in~\cite{CaRi19}, wherein a machine-learned regression model was used to learn a direct relationship between the internal state of an element and its forces. The approach enforces frame indifference and conservation of linear and angular momentum and was shown to reduce the computational cost as opposed to a traditional finite element model. Refs.~\cite{BaMaSt22,TaMeSt22} proposed additional techniques for computing the tangent stiffness matrix of smart elements based on finite differences, automatic differentiation, and neural networks. In addition,~\cite{TaMeSt22} proposes a composite loss function with terms for both the element forces as well as the tangent stiffness. Similarly,~\cite{KoBaMa20} proposes nonlinear meta elements (i.e., element patches) using deep learning. Unlike~\cite{CaRi19,TaMeSt22} this approach learns generalizable meta elements (where each meta element comprises multiple finite elements). The primary input to the regression model defining the meta element is the displacement field on element boundary, while the outputs are the displacements, stresses, and forces within the meta element. Lastly,~\cite{BeZaRo23} explores quantifying uncertainties in machine-learned elements that map from displacements to forces by virtue of an extension to deep ensembles~\cite{LaBePr17}. 

The aforementioned approaches seek to learn a generalizable (meta-)element that can be applied repeatedly throughout the domain. Another similar body of work focuses on \textit{domain-decomposition} in where a computational domain is divided into subdomains. Once partitioned, each subdomain can be solved independently, and a global solution is obtained by coupling the various domains, e.g., with the Schwarz method~\cite{MoTeAl17}. 
Data-driven techniques have received significant attention in this field as well, where the general idea is to use a cheaper-to-evaluate surrogate to model a complex subdomain. Refs.~\cite{LiTaWu20,WuXuYi20}, for instance, proposed D3M and DeepDDM, both of which employ physics-informed neural networks (PINNs) to learn subdomain-solutions conditioned on boundary data, and demonstrated their approaches on Poisson's equation and Schrodinger's equation. We note that D3M and DeepDDM employ PINNs for all subdomains, as opposed to a FEM for some of the subdomains. Refs.~\cite{MaRoEi02,HuKnPa13,IaQuRo16,HuCa21,CaKuTe22,BaTeMo22,IoSaTa23,CaBoKu23} propose domain decomposition strategies leveraging reduced basis and proper-orthogonal decomposition (Petrov-)Galerkin reduced-order models. These methods operate by projecting the finite element model onto low-dimensional subspaces computed, e.g., from training data. The different subdomains are then coupled via, e.g., Lagrange multipliers~\cite{HuCa21}, the Schwarz method~\cite{MoTeAl17,MoTePh22}. Lastly, one final approach of relevance was proposed in~\cite{AhSaKa20} for multiscale problems, where the authors aim to resolve only the domain of interest by learning an interface condition leveraging long short term memory (LSTM) networks and ideas from upwinding. The setup was demonstrated on the transient Burgers' equation and Euler equations in one dimension.

In the present work, we develop a machine learning strategy for replacing portions of a computational domain, such as threaded fasteners, with a data-driven surrogate.  We develop the proposed strategy within the context of hypoelastic material models for static solid mechanics (SM) FEM simulations of bodies under finite deformations exposed to contact and friction; we leave the extension of the approach to nonlinear material models and dynamic simulations as future work. In the proposed strategy we decompose a computational domain into an ``outer" domain that we resolve with an FEM, and an ``inner" domain. We then develop a machine-learned model for the impact of the inner domain on the outer domain. In essence, for SM our machine-learned surrogate learns the map from displacements on the interface boundary to forces contributed by the inner domain to the outer domain on the same interface boundary (we refer to these forces as internal forces). 

We propose two displacement-to-force mappings in the present work. First, we propose a direct displacement-to-force mapping with dimension reduction. In this approach, we identify low-dimensional subspaces for compact approximations to the interface displacements and internal forces via proper orthogonal decomposition (POD)~\cite{BeHoLu93}. A mapping between reduced coordinates for these fields is then learned via regression; i.e., for interface displacements $\mathbf{u}$ and internal forces $\mathbf{f}$, we learn a mapping $\mathbf{f} = \mathbf{g}(\mathbf{u})$ where $\mathbf{g}$ is the learned surrogate. The second approach we propose is a structure-preserving displacement-to-stiffness mapping with dimension reduction. In this approach, we again identify low-dimensional structure for the displacements and internal forces on the interface, but instead learn a displacement-to-force mapping in the form of a stiffness matrix, i.e., for interface displacements $\mathbf{u}$ and internal forces $\mathbf{f}$, we learn a mapping $\mathbf{f} = \mathbf{K}(\mathbf{u}) \mathbf{u}$ where $\mathbf{K}(\mathbf{u})$ is a learned stiffness which can be a nonlinear function of the displacements. We further enforce the learned stiffness matrix to be symmetric positive semi-definite (SPSD)\footnote{We note that our ML model is only learning a boundary coupling term and enforcing SPSD structure is sufficient to guarantee that the full assembled stiffness matrix in the coarse-scale problem is symmetric positive definite.}. We show, in a simplified setting, that the addition of this SPSD model results in a coarse-scale problem that is symmetric positive definite (SPD) and well-posed. Further, enforcing SPSD structure allows our ML model to be interfaced with common SM solvers such as conjugate gradient. Thematically similar ideas were pursued in Ref.~\cite{XuHuDa21,AsAvFa22} within the context of learning constitutive relations. As will be shown, we find this latter approach is critical for successful integration of the method into our finite element solver.

The current work has commonalities with several existing efforts. For example, the present approach can be viewed as a smart element or a ``meta-element" strategy in which we are learning a meta-element for our inner domain(s). This synergy is particularly relevant if the inner domain is a component which is repeated throughout the computational domain, like a threaded fastener. We propose a more limited approach than~\cite{CaRi19,TaMeSt22,KoBaMa20}, however: we aim to remove certain hard-to-resolve portions of a domain rather than develop a general smart element. In this sense, our proposed methodology can be interpreted as a type of Guyan reduction that is applicable to the nonlinear regime. Second, our approach can be viewed as a domain-decomposition method, where we employ an ML model to simulate the impact of fine-scale domain on the coarse-scale domain. Comparing the present work to the literature, we highlight several novelties:
\begin{enumerate}
\item  We propose regression models equipped with a dimension reduction strategy to reduce both the training cost of the model as well as the number of datum required to learn an accurate and generalizable model. This strategy enables our approach to be applied to systems where the ML model must predict thousands of degrees of freedom.
\item We propose two types of displacement-to-force mappings: one that directly maps from displacement to forces, and one that maps from displacements to forces by virtue of learning an SPSD stiffness matrix. This second approach is novel. Refs.~\cite{CaRi19,TaMeSt22,KoBaMa20}, for example, only consider direct displacement-to-force mappings, and the stiffness matrix associated with these mappings is extracted via, e.g., automatic differentiation. By directly learning a stiffness matrix we can enforce structure in our learned models. Our numerical analysis and experiments demonstrate that this embedded structure results in increased robustness. Additionally, the proposed displacement-to-stiffness mapping addresses the issue of stiffness computation studied in~\cite{TaMeSt22}. In our setting the stiffness is a direct output of the neural network and does not need to be probed with finite differences, computed via back propagation, or estimated with an auxiliary network.
\end{enumerate}  
This manuscript proceeds as follows. In Section~\ref{sec:formulation} we provide the general formulation for the governing equations of solid mechanics, while in Section~\ref{sec:multiscale} we present an equivalent decomposed formulation. Section~\ref{sec:ml-framework} outlines our machine-learning framework, and Section~\ref{sec:simple_analysis} provides analysis of the coupled FEM-ML systems when applied to the equations of linear elasticity in one dimension. Section~\ref{sec:experiment} provides numerical experiments, while Section~\ref{sec:conclude} provides conclusions and perspectives.



\section{Monolithic formulation}\label{sec:formulation}
We begin by outlining the governing equations for finite-deformation solid mechanics. While we consider physics such as friction and contact in our numerical experiments, these details are not presented here for brevity. We refer the interested reader to the Sierra/SM theory manual~\cite{sierra} for more details.

Consider a body defined on the reference domain $\domain \subset \RR{3}$ with boundary $\domainBoundary$ undergoing finite deformations described by $\eulerianCoordinate = \deformationMap(\materialCoordinate)$, where $\deformationMap : \domain \rightarrow \deformedDomain$ is the deformation map, $\materialCoordinate \in \domain$ are the reference coordinates, $\eulerianCoordinate \in \deformedDomain$ are the deformed coordinates, and $\deformedDomain = \deformationMap(\domain)$ is the deformed domain. 
In what follows we use $\deformationMap \equiv \deformationMap\left(\materialCoordinate \right)$ for brevity. The governing equations for a quasi-static mechanical problem are given as 
\begin{equation}\label{eq:forceBalanceNL}
\begin{split}
 -\nabla \cdot \PKStress \left(\deformationMap \right) & = \forcing \text{ on } \domain. 
\end{split}
\end{equation}
In the above, $\PKStress : \deformationMap \mapsto J(\deformationMap) \CauchyStress(\deformationMap) \deformationGradient^{-T}(\deformationMap)$ is the first Piola–Kirchhoff stress tensor, $\deformationGradient(\deformationMap) = \nabla \deformationMap$ is the deformation gradient, $J(\deformationMap) = \mathrm{det} \left(\deformationGradient\left(\deformationMap\right) \right)$, $\CauchyStress(\deformationMap)$ is the Cauchy stress tensor, and $\forcing : \domain \rightarrow  \RR{3}$ is the body forcing per unit reference volume. We emphasize that the first Piola-Kirchoff stress tensor is a nonlinear function of the deformation map and contains the constitutive relationship. We assume the domain boundary $\domainBoundary$ can be decomposed as $\domainBoundary  = \domainBoundaryDirichlet\cup \domainBoundaryTraction$ with $\domainBoundaryDirichlet \cap \domainBoundaryTraction = \emptyset$, where $\domainBoundaryDirichlet$ comprises a Dirichlet boundary and $\domainBoundaryTraction$ a traction boundary such that Eq.~\eqref{eq:forceBalanceNL} is supplemented with the boundary conditions
\begin{equation*}
\begin{split}
\PKStress(\deformationMap) \cdot \normal &=  \tractionBC \text{ on } \domainBoundaryTraction \\
\deformationMap & = \deformationMapBC \text{ on } \domainBoundaryDirichlet.
\end{split}
\end{equation*}
In the above, $\tractionBC$ denotes the prescribed boundary tractions, $\deformationMapBC$ the prescribed Dirichlet boundary conditions, and $\normal$ defines the normal vector to the boundary $\domain$.  

We consider the standard weighted residual form of Eq.~\eqref{eq:forceBalanceNL}, which reads as follows: find $\deformationMap \in \trialSpace$ such that
\begin{equation}\label{eq:forceBalanceVariationalNL}
\int_{\domain} \PKStress \left(\deformationMap\right) : \nabla \testFunction dV  - \int_{\domainBoundaryTraction} (\PKStress\left(\deformationMap \right) \cdot \normal) \cdot \testFunction dS    = \int_{\domain} \forcing \cdot \testFunction dV , \qquad \forall \testFunction \in \testSpace,
\end{equation}
where $\trialSpace$ denotes a suitable vector-valued trial space on $\domain$ satisfying the Dirichlet boundary condition $\deformationMapBC$ on $\domainBoundaryDirichlet$ and $\testSpace$ denotes a suitable vector-valued test space on $\domain$ vanishing on $\domainBoundaryDirichlet$.
The Galerkin method with Lagrange finite elements comprises a standard approach to solve Eq.~\eqref{eq:forceBalanceVariationalNL}. To this end, let $\trialSpaceFEM \subset \trialSpace$ denote a conforming trial space obtained via a Lagrange finite element discretization of the reference domain $\referenceDomain$ into $\nElements$ non-overlapping elements $\domain_k, k=1,\ldots,\nElements$, and let $\testSpaceFEM$ be the corresponding test space. The spatially discrete counterpart to Eq.~\eqref{eq:forceBalanceVariationalNL} reads: find $\deformationMapFEM \in \trialSpaceFEM$ such that
\begin{equation}\label{eq:forceBalanceGalerkinNL}
\int_{\referenceDomain} \PKStress \left( \deformationMapFEM\right)  : \nabla \testFunctionFEM dV  - \int_{\domainBoundaryTraction} \left( \PKStress \left(\deformationMapFEM\right) \cdot \normal \right) \cdot \testFunctionFEM dS    = \int_{\referenceDomain} \forcing \cdot \testFunctionFEM dV , \qquad \forall \testFunctionFEM \in \testSpaceFEM.
\end{equation}

Equation~\eqref{eq:forceBalanceGalerkinNL} comprises a monolithic set of governing equations. It is often the case, however, that the body $\domain$ is made up of multiple subdomains, where each subdomain may require, e.g., different resolution requirements. 
One such example that motivates the present work is systems-level FEM models comprising threaded fasteners. In this case, the threaded fastener model can require a much finer discretization than the remainder of the domain and can be a computational bottleneck. 
To mitigate this issue we propose a data-driven modeling approach.

\section{Domain decomposition formulation}\label{sec:multiscale} 
\begin{figure}
\begin{center}
\begin{subfigure}[t]{0.65\textwidth}
\includegraphics[trim={7cm 7cm 7cm 7cm},clip,width=1.0\linewidth]{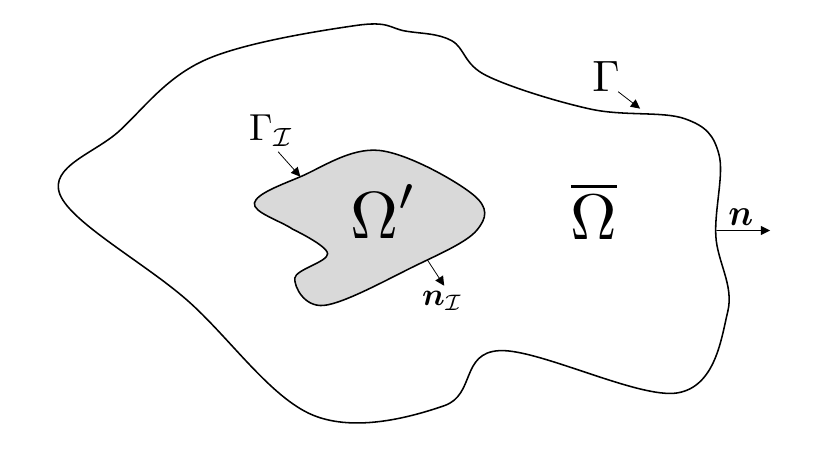}
\end{subfigure}
\caption{Simplified depiction of domain-decomposition formulation.}
\label{fig:domainDecomposition}
\end{center}
\end{figure}

We consider a decomposition of the reference domain $\referenceDomain$ into two non-overlapping components such that $\referenceDomain = \domainOne \cup \domainTwo$ with the interface $\domainBoundaryOneTwo = \domainOne \cap \domainTwo$ as in Figure~\ref{fig:domainDecomposition}. We take $\domainOne$ to be an ``outer" coarse-scale domain which we aim to solve directly with a finite element method, while we take $\domainTwo$ to be an ``inner" fine-scale domain. We aim to not resolve $\domainTwo$, but rather model its impact on $\domainOne$. To this end, we introduce a decomposition of the deformation map on $\domainOne$ and $\domainTwo$ 
$$\deformationMap = \deformationMapOne + \deformationMapTwo$$
such that $\deformationMapOne = \bz$ on $\domainTwo$, $\deformationMapTwo = \bz$ on $\domainOne$, and $\deformationMapOne = \deformationMapTwo \equiv \deformationMap$ on $\domainBoundaryOneTwo$. 
We note that, by definition, $\PKStress(\deformationMapOne + \deformationMapTwo) = \PKStress(\deformationMapOne)$ on $\domainOne$ and $\PKStress(\deformationMapOne + \deformationMapTwo) = \PKStress(\deformationMapTwo)$ on $\domainTwo$.
With this decomposition, a coarse-scale problem on $\domainOne$ can be written as: find $\deformationMapOne \in \trialSpaceOne$ such that $\forall \testFunctionOne \in \testSpaceOne$
\begin{equation}\label{eq:coarseEquation}
 \int_{\domainOne} \PKStress\left(\deformationMapOne \right): \nabla \testFunctionOne dV  - \int_{\domainBoundaryTraction} \left(\PKStress\left(\deformationMapOne\right) \cdot \normal \right) \cdot \testFunctionOne dS   = \int_{\domainOne} \forcing \cdot \testFunctionOne dV  - \int_{\domainBoundaryOneTwo}   \left( \PKStress\left(\deformationMapTwo \right)\cdot \normalOneTwo \right)  \cdot \testFunctionOne  dS ,
\end{equation}
where again $\trialSpaceOne$ and $\testSpaceOne$ are suitable vector-valued trial and test spaces on $\domainOne$ and $\normalOneTwo$ is the outward normal vector on the interface $\domainBoundaryOneTwo$ (see Figure~\ref{fig:domainDecomposition}). Our goal is to solve the coarse-scale equation. Let $\coarseTrialSpaceFEM \subset \trialSpaceOne$ denote a conforming trial space obtained via a Lagrange finite element discretization of the domain $\domainOne$ into $\nCoarseElements$ non-overlapping elements $\domainOne_k, k=1,\ldots,\nCoarseElements$ with the node set $\coarseNodeSet$, $| \coarseNodeSet| = \nCoarseNodes$, and let $\coarseTestSpaceFEM$ denote the corresponding test space. For notational purposes, we additionally denote the interface boundary node set as $\coarseNodeSetDomainBoundaryOneTwo \subset \coarseNodeSet$ with $| \coarseNodeSetDomainBoundaryOneTwo | = \nCoarseNodesDomainBoundaryOneTwo$.
The spatially discrete counterpart to Eq.~\eqref{eq:coarseEquation} reads: find $\deformationMapCoarseFEM \in \coarseTrialSpaceFEM$ such that $\forall \testFunctionOneFEM \in \coarseTestSpaceFEM$
\begin{equation}\label{eq:coarseEquationFEM}
\begin{split}
 \int_{\domainOne} \PKStress\left(\deformationMapCoarseFEM\right) : \nabla \testFunctionOneFEM dV  - \int_{\domainBoundaryTraction} \left(\PKStress\left(\deformationMapCoarseFEM \right)\cdot \normal\right) \cdot \testFunctionOneFEM dS   = \\ \int_{\domainOne} \forcing \cdot \testFunctionOneFEM dV  - \int_{\domainBoundaryOneTwo}   \left( \PKStress\left(\deformationMapTwo \right) \cdot \normalOneTwo \right)  \cdot \testFunctionOneFEM  dS.
\end{split} 
\end{equation}
The above can be written as a set of $\nCoarseDofs =  3 \nCoarseNodes$ algebraic equations by introducing the vector-valued trial basis functions $\trialFunctionOneFEM_i$, $i=1,\ldots,\nCoarseNodes$ and vector-valued test basis functions ${\testFunctionOneFEM}_i$, $i=1,\ldots,\nCoarseNodes$ such that $\mathrm{span}\{ {\trialFunctionOneFEM}_i \}_{i=1}^{\nCoarseNodes} = \coarseTrialSpaceFEM$, $\mathrm{span}\{ {\testFunctionOneFEM}_i \}_{i=1}^{\nCoarseNodes} = \coarseTestSpaceFEM$. 
The system~\eqref{eq:coarseEquationFEM} simplifies to a force balance at each node in the case of Lagrange finite elements. 
We observe that the system~\eqref{eq:coarseEquationFEM} is unclosed as it requires specification of the boundary term $  \left( \PKStress\left(\deformationMapTwo \right) \cdot \normalOneTwo \right)  \cdot \testFunctionOneFEM  $, i.e., we need to specify the forcing exhibited by the inner domain on the outer domain. We employ a machine-learning approach for this purpose.

\section{Machine-learning framework}\label{sec:ml-framework}
Developing a method that only requires solving the coarse-scale system requires a model for the impact of the fine scales on the coarse scales. To this end, we propose a coupled FEM-ML formulation wherein the ML model is used to describe the impact of the fine scales on the coarse scales at the interface boundary. This section outlines our coupled FEM-ML formulation, describes the offline-online approach employed to construct the ML model, and lastly outlines several candidate machine-learning models.  
\subsection{FEM-ML coupled formulation}
Closing the coarse-scale equation requires specification of the traction applied by the fine scales on the domain boundary, i.e., we require a model for the term $\int_{\domainBoundaryOneTwo}   \left( \PKStressFine \cdot \normalOneTwo \right)  \cdot \testFunctionOneFEM  dS.$ We propose a machine learning approach for this purpose, i.e., we propose to develop a machine learning model $\mlModelVec(\cdot)  = \left[ \mlModel_1(\cdot)^T ,\ldots,\mlModel_{\nCoarseNodes}\left(\cdot \right)^T \right]^T \in \RR{\nCoarseDofs}$ such that 
$$\mlModel_i(\features ) \approx  \int_{\domainBoundaryOneTwo}   \left( \PKStress\left(\deformationMapTwo\right) \cdot \normalOneTwo \right)  \cdot {\testFunctionOneFEM}_i dS \in \RR{3},$$
for $i=1,\ldots,\nCoarseNodes$. In the above, $\mlModel_i(\features)$ denotes the ML model output for the $i$th test function for the $3$ DOFs (i.e., DOFs in the $x_1,x_2$, and $x_3$ directions) and $\features$ are the to-be-defined features fed into the ML model.
The coarse-scale equation then becomes: find $\deformationMapCoarseFEM \in \coarseTrialSpaceFEM$ such that for $i=1,\ldots,\nCoarseNodes$,
\begin{equation}\label{eq:coarseEquationFEMML}
 \int_{\domainOne} \PKStress\left(\deformationMapCoarseFEM\right) : \nabla {\testFunctionOneFEM}_i dV  - \int_{\domainBoundaryTraction} \left(\PKStress\left(\deformationMapCoarseFEM \right) \cdot \normal\right) \cdot {\testFunctionOneFEM}_i dS   + \mlModel_i(\features)= \int_{\domainOne} \forcing \cdot {\testFunctionOneFEM}_i dV . 
\end{equation}


\subsection{Offline-online workflow}
To develop the data-driven model, we employ an \textit{offline--online} workflow. In the offline phase, we perform a computationally intensive process that involves solving the \textit{full} finite element model~\eqref{eq:forceBalanceGalerkinNL} for $\nSamples$ different configurations (e.g., different boundary conditions) to obtain $\nSamples$ deformation maps $\{\deformationMapFEM^i\}_{i=1}^{\nSamples}$, where $\nSamples$ are the number of training samples.
These sample deformation maps are then processed to compute the target training matrix $\forceSample \in \RR{ \nCoarseDofs \times \nSamples}$ comprising the terms $\int_{\domainBoundaryOneTwo} \left({\PKStress \left(\deformationMapFEM^j \right)} \cdot \normal \right) \cdot {\testFunctionOneFEM}_i dS$ for $i=1,\ldots,\nCoarseNodes$, $j=1,\ldots,\nSamples$. 
Analogously, we collect the feature training matrix $\featuresMatrix \in \RR{\nFeatures \times \nSamples}$, where $\nFeatures$ are the number of features. 
In this work we consider the case where the features are the nodal values of the displacement fields, i.e., 
$$\features = \displacementCoarseDiscreteFEM \equiv  \left[ \left[\deformationMapCoarseFEM(\materialCoordinate^{\domainOne}_1) - \materialCoordinate_1^{\domainOne} \right]^T \; \cdots \; \left[\deformationMapCoarseFEM(\materialCoordinate^{\domainOne}_{\nCoarseNodes}) - \materialCoordinate_{\nCoarseNodes}^{\domainOne} \right]^T \right]^T \in \RR{\nCoarseDofs},$$
where $\materialCoordinate^{\domainOne}_{i}$, $i=1,\ldots,\nCoarseNodes$ denotes coordinates of the $i$th coarse node. 

\subsubsection{Simplification for Lagrange finite elements}
In the present work we exclusively consider application to Lagrange finite elements. In this case the term $\int_{\domainBoundaryOneTwo} \left({\PKStress}\left(\deformationMapTwo \right) \cdot \normal \right) \cdot {\testFunctionOneFEM}_i dS$ is exactly zero for basis functions not on the interface $\domainBoundaryOneTwo$ and we only need to construct a model for basis functions on the interface such that $\mlModelVec(\cdot) \in \RR{\nCoarseDofsDomainBoundaryOneTwo}$, where $\nCoarseDofsDomainBoundaryOneTwo = 3 \nCoarseNodesDomainBoundaryOneTwo$. In this case we need only collect the target training matrix for the $\nCoarseNodesDomainBoundaryOneTwo$ basis functions on the interface. This results in a target training matrix of size $\forceSample \in \RR{ \nCoarseDofsDomainBoundaryOneTwo \times \nSamples}$. Analogously, for features we employ only the nodal values of the displacement fields on the interface, i.e., $\features = \displacementCoarseDiscreteFEM^{\domainBoundaryOneTwo}$ such that $\featuresMatrix \in \RR{\nCoarseDofsDomainBoundaryOneTwo \times \nSamples}$. This approach will be adapted for the remainder of the manuscript. 


\subsection{Machine-learning model structure}
Having defined the feature matrix $\featuresMatrix$ and the response matrix $\forceSample$ we now outline two candidate approaches for the ML model $\mlModelVec$.
\vspace{0.2 in}\\
\noindent \textbf{Approach 1: Direct mapping with proper orthogonal decomposition.} 
The first approach we consider comprises learning a direct feature-to-response mapping with the addition of POD to reduce the dimensionality of the input feature space and the response space. This approach can be described as
$$\mlModelVec : \displacementCoarseDiscreteFEM^{\domainBoundaryOneTwo} \mapsto  \basisForces \reducedMlModelVec \left( \basisDisplacements^T \displacementCoarseDiscreteFEM^{\domainBoundaryOneTwo} ; \weights \right)  + \forceOffset,$$
where $\forceOffset \in \RR{\nCoarseDofsDomainBoundaryOneTwo}$ is a constant vector (accounting, e.g., for body forces, preload), $\basisForces \in \RR{\nCoarseDofsDomainBoundaryOneTwo \times \reducedDimensionForces}$ and $\basisDisplacements \in \RR{\nCoarseDofsDomainBoundaryOneTwo \times \reducedDimensionDisplacement}$ are \textit{basis} matrices for the interface forces and displacements, respectively, obtained through POD, $\reducedDimensionDisplacement$ and $\reducedDimensionForces$ are the reduced dimension for the displacements and forces, both of which are model hyperparameters,  $\reducedMlModelVec : \RR{\reducedDimensionDisplacement} \times \RR{\nWeights} \rightarrow \RR{\reducedDimensionForces}$ comprises a machine-learned model that maps from reduced displacement states to reduced force states, and $\weights \in \RR{\nWeights}$ are model parameters (obtained from training). The basis matrices for forces and displacements are obtained by solving the optimization problems
\begin{equation}\label{eq:basis_force}
\basisForces = \underset{ \basisForces^* \in \{ \RR{\nCoarseDofsDomainBoundaryOneTwo \times \reducedDimensionForces} | [\basisForces^*]^T \basisForces^* = \mathbf{I} \} }{\text{arg min}} \| \basisForces^* [\basisForces^*]^T \forceSampleOffset - \forceSampleOffset \|_2^2,
\end{equation}
\begin{equation}\label{eq:basis_displacement}
\basisDisplacements = \underset{ \basisDisplacements^* \in \{ \RR{\nCoarseDofsDomainBoundaryOneTwo \times \reducedDimensionDisplacement} | [\basisDisplacements^*]^T \basisDisplacements^* = \mathbf{I} \} }{\text{arg min}} \| \basisDisplacements^* [\basisDisplacements^*]^T \featuresMatrix - \featuresMatrix \|_2^2,
\end{equation}
where $\forceSampleOffset = \forceSample - \left[ \forceOffset \; \cdots \; \forceOffset \right].$
These optimization problems can be solved efficiently with a (thin) singular value decomposition (SVD)~\cite{Vo11}. The algorithm is described in Appendix~\ref{appendix:pod}. We emphasize that by employing POD, the training of the model scales with the reduced dimensions $\reducedDimensionDisplacement$ and $\reducedDimensionForces$ rather than the number of interface nodes.


A variety of regression models can be used to learn the reduced model $\reducedMlModelVec$. Here we investigate two approaches spanning different regimes of model complexity: linear regression and neural networks. We investigate linear regression due to its simplicity and interpretability, while we investigate neural networks due to their high capacity to learn complex functions. These two strategies are now described.
\begin{itemize}
\item \textbf{Linear least squares (LLS)}. In a linear regression approach the reduced ML model is defined by a matrix $\reducedLinearRegressionMatrix \in \RR{\reducedDimensionForces \times \reducedDimensionDisplacement}$ such that the full model can be written as
$$\mlModelVec : \displacementCoarseDiscreteFEM^{\domainBoundaryOneTwo} \mapsto  \basisForces \reducedLinearRegressionMatrix \basisDisplacements^T \displacementCoarseDiscreteFEM^{\domainBoundaryOneTwo}.$$
The reduced system matrix is defined by the solution to the linear least-squares problem
\begin{equation}\label{eq:reducedLeastSquaresMatrix}
\reducedLinearRegressionMatrix = \underset{ \reducedLinearRegressionMatrix^* \in \RR{\reducedDimensionForces \times \reducedDimensionDisplacement  } }{\text{arg min}} \| \basisForces^T \forceSampleOffset - \reducedLinearRegressionMatrix^* \basisDisplacements^T \featuresMatrix \|_2^2.
\end{equation}
The least-squares problem can be efficiently solved as it can be broken down into $\reducedDimensionForces$ least-squares problems for each row of $\reducedLinearRegressionMatrix$. 

\item \textbf{Neural network (NN)}
This method learns a feed forward neural network for the mapping from the (reduced) displacement field to the (reduced) force field. In this case the ML model is given by
$$\mlModelVec : \displacementCoarseDiscreteFEM^{\domainBoundaryOneTwo} \mapsto  \basisForces \reducedNeuralNetwork \left( \basisDisplacements^T \displacementCoarseDiscreteFEM^{\domainBoundaryOneTwo} ; \weights \right) $$
where $\reducedNeuralNetwork : \RR{\reducedDimensionDisplacement} \times \RR{\nWeights} \rightarrow \RR{\reducedDimensionForces}$ comprises a neural network model that maps from the reduced displacement states to the reduced force states and $\weights \in \RR{\nWeights}$ are the network weights. In the present work, these network weights are obtained from approximate minimization of the loss function,
$$\underset{ \weights \in \RR{\nWeights}}{\mathrm{minimize}} \sum_{i=1}^{\nSamples} \|   \reducedMlModelVec\left( \basisDisplacements^T \featuresMatrix_i ; \weights\right) - \basisForces^T \forceSampleOffset_i \|_2^2.$$ 

\end{itemize}

\noindent \textbf{Approach 2: Stiffness matrix mapping with proper orthogonal decomposition.} 
The approach outlined above is the simplest learning framework one can employ. This approach, however, does not directly enforce any physical or mathematical properties of the underlying FEM. As will be shown later in this manuscript, we found that embedding the above models into an FEM code often led to unstable solutions. To address  this challenge, we propose a second approach that embeds structure into the model by learning an SPSD stiffness matrix. We again use POD to reduce the dimensionality of the input feature and response spaces. This approach can be described as
\begin{equation}\label{eq:stiffnessApproach}
\mlModelVec : \displacementCoarseDiscreteFEM^{\domainBoundaryOneTwo} \mapsto  \basisCombined \left[ \reducedMlModelStiffnessVec \left( \basisCombined^T \displacementCoarseDiscreteFEM^{\domainBoundaryOneTwo} ; \weights \right) \right]  \left[ \reducedMlModelStiffnessVec \left( \basisCombined^T \displacementCoarseDiscreteFEM^{\domainBoundaryOneTwo} ; \weights \right) \right]^T  \basisCombined^T \displacementCoarseDiscreteFEM^{\domainBoundaryOneTwo} + \forceOffset 
\end{equation} 
where $\reducedMlModelStiffnessVec : \RR{\reducedDimensionCombined} \times \RR{\nWeights} \rightarrow \RR{\reducedDimensionCombined \times \reducedDimensionCombined}$ comprises a machine-learned model that maps from the reduced displacement states to the lower triangular part of a reduced stiffness matrix, and $\basisCombined \in \RR{\nCoarseDofsDomainBoundaryOneTwo \times \left( \reducedDimensionForces + \reducedDimensionDisplacement\right) } = \textit{orthogonalize}\left( \left[ \basisForces , \basisDisplacements \right] \right)$ is a basis matrix combining the (orthogonalized) union of the force basis and displacement basis. We note the orthogonalization can be performed efficiently with a QR decomposition. 

Critically, Eq.~\eqref{eq:stiffnessApproach} can be written in matrix form as 
$$\mlModelVec \left(  \displacementCoarseDiscreteFEM^{\domainBoundaryOneTwo}  \right) = \stiffnessMatrix_{\mathrm{ML}} \left(   \displacementCoarseDiscreteFEM^{\domainBoundaryOneTwo}  \right) \displacementCoarseDiscreteFEM^{\domainBoundaryOneTwo}  + \forceOffset, $$
where $ \stiffnessMatrix_{\mathrm{ML}} (   \displacementCoarseDiscreteFEM^{\domainBoundaryOneTwo}  ) =  [\stiffnessMatrix_{\mathrm{ML}} (\displacementCoarseDiscreteFEM^{\domainBoundaryOneTwo})]^T \in \RR{\nCoarseDofsDomainBoundaryOneTwo \times \nCoarseDofsDomainBoundaryOneTwo}$ is an SPSD stiffness matrix. This structure is directly enforced in the model by learning the lower triangular matrix $\reducedMlModelStiffness$. Enforcing our model to be SPSD makes the coupled FEM-ML system more amenable for conjugate gradient solvers, which are commonly employed in solid mechanics. Lastly, we note the use of the combined basis $\basisCombined$ for both the force and displacement field, opposed to the individual bases employed in Approach 1. The reason for employing the combined basis is that both the force field and displacement field must employ the same basis for the model~\eqref{eq:stiffnessApproach} to be SPSD\footnote{For example, if we employ different bases for the force and displacement, then stiffness matrix associated with~\eqref{eq:stiffnessApproach} would take the form 
$
 \basisForces \left[ \reducedMlModelStiffnessVec \left( \basisDisplacements^T \displacementCoarseDiscreteFEM^{\domainBoundaryOneTwo} ; \weights \right) \right]  \left[ \reducedMlModelStiffnessVec \left( \basisDisplacements^T \displacementCoarseDiscreteFEM^{\domainBoundaryOneTwo} ; \weights \right) \right]^T  \basisDisplacements^T 
$
which is not symmetric.}
. In the present work we achieve this by combining both the bases for the force field and displacement fields. For simplicity we consider the case where $\reducedDimensionCombined = \reducedDimensionDisplacement + \reducedDimensionForces$; it is possible that a more compact representation can be obtained by considering a single POD problem where the displacements and forces are stacked. 

A variety of regression models can again be used to learn the reduced stiffness matrix $\reducedMlModelStiffness$. We again consider a linear regression approach and a neural-network-based approach.
\begin{itemize}
\item \textbf{Symmetric positive semi-definite linear least-squares (SPSD-LLS)}. This approach learns a linear model for the stiffness matrix. This method is described by
$$\mlModelVec : \displacementCoarseDiscreteFEM^{\domainBoundaryOneTwo} \mapsto  \basisCombined \reducedLowerLinearRegressionMatrix \reducedLowerLinearRegressionMatrix^T \basisCombined^T \displacementCoarseDiscreteFEM^{\domainBoundaryOneTwo} + \forceOffset$$
where $\reducedLowerLinearRegressionMatrix \in \RR{\reducedDimensionCombined \times \reducedDimensionCombined}$ is defined by the optimization problem
\begin{equation}\label{eq:leastSquaresMatrixLower}
\reducedLowerLinearRegressionMatrix = \underset{ \reducedLowerLinearRegressionMatrix^* \in \RR{\reducedDimensionCombined \times \reducedDimensionCombined } }{\text{arg min}} \| \basisCombined^T \forceSampleOffset - \reducedLowerLinearRegressionMatrix^*  \left[ \reducedLowerLinearRegressionMatrix^*\right]^T \basisCombined^T \featuresMatrix \|_2^2.
\end{equation}
By definition, SPSD-LLS learns a system matrix that is SPSD. To the best of our knowledge, no analytic solution exists to the optimization problem~\eqref{eq:leastSquaresMatrixLower}, but it can be solved in a straightforward manner by casting it as a constrained optimization problem. In the present work we solve the optimization problem~\eqref{eq:leastSquaresMatrixLower} with the optimization package \textsf{cvxpy}~\cite{diamond2016cvxpy,agrawal2018rewriting}. 

\item \textbf{Symmetric positive semi-definite neural network regression (SPSD-NN)}
This approach learns a feed forward neural network for the mapping from the (reduced) displacement field to the (reduced) lower triangular stiffness matrix. In this case the ML model is given by
$$\mlModelVec : \displacementCoarseDiscreteFEM^{\domainBoundaryOneTwo} \mapsto  \basisCombined \left[ \reducedNeuralNetworkStiffness \left( \basisCombined^T \displacementCoarseDiscreteFEM^{\domainBoundaryOneTwo} ; \weights \right) \right]  \left[ \reducedNeuralNetworkStiffness \left( \basisCombined^T \displacementCoarseDiscreteFEM^{\domainBoundaryOneTwo} ; \weights \right) \right]^T  \basisCombined^T \displacementCoarseDiscreteFEM^{\domainBoundaryOneTwo}   + \forceOffset $$
where $\reducedNeuralNetworkStiffness : \RR{\reducedDimensionCombined} \times \RR{\nWeights} \rightarrow \RR{\reducedDimensionCombined \times \reducedDimensionCombined}$ is  the neural network model for the reduced lower stiffness matrix.   
The network weights are obtained from approximate minimization of the loss function,
$$\underset{ \weights \in \RR{\nWeights}}{\mathrm{minimize}} \sum_{i=1}^{\nSamples} \| 
\left[ \reducedNeuralNetworkStiffness \left( \basisCombined^T  \featuresMatrix_i ; \weights \right) \right]  \left[ \reducedNeuralNetworkStiffness \left( \basisCombined^T  \featuresMatrix_i ; \weights \right) \right]^T  \basisCombined^T  \featuresMatrix_i - \basisCombined^T \forceSampleOffset_i \|_2^2.$$ 
\end{itemize}

In summary, we considered four model forms: a linear model (LLS) and nonlinear (NN) model that directly map from the interface displacements to the interface forces, as well as a linear model (SPSD-LLS) and nonlinear model (SPSD-NN) that directly map from interface displacements to an SPSD stiffness matrix. 

\begin{remark}
We emphasize that the above approaches learn the mapping between displacements and forces with the purpose of removing a subdomain from the computational domain. A similar approach could be to learn, e.g., the strain--stress mapping on the domain boundary. This type of approach was employed in Refs.~\cite{XuHuDa21,AsAvFa22}, which leverage neural networks to learn the constitutive relationship defining the stress-strain relationship at every point a computational domain. In the present work we directly pursue a displacement--force mapping for ease of implementation into our finite element solver. Specifically, the contribution of our ML model simply appears as an additional (state-dependent) source term in the residual of the governing equations (as opposed to, e.g., a stress field that needs to be integrated). 
\end{remark}

\begin{remark}
The above formulations are designed for path-independent static problems and will be incomplete for dynamic path-dependent problems (e.g., plasticity) where history effects can be important.
\end{remark} 

\section{Positive definiteness of coarse-scale model in a simplified setting}\label{sec:simple_analysis}
In this section we demonstrate that, in the simplified setting of linear elasticity in one-dimension with homogeneous boundary conditions, the SPSD stiffness-based models result in a global coarse scale system that is SPD and, as a result, have invertible system matrices with unique solutions. We further show that, for the direct displacement-to-force model, we cannot guarantee an invertible system with a unique solution. 

The governing equations of linear elasticity in one dimension are
\begin{equation}\label{eq:oneDimensionalExample}
-A E \frac{\partial^2 \displacementScalar}{\partial x^2} = \forcingScalar
\end{equation}
for $x \in [0,L]$ where $\displacementScalar : [0,L] \rightarrow \RR{}$ is the displacement with boundary conditions $\displacementScalar(0) = \displacementScalar(L) = 0$. In the above, $A \in \RR{^+}$ is the cross-sectional area, $E \in \RR{+}$ is the modulus of elasticity, and $\forcingScalar \in \RR{}$ is a constant body forcing term. Discretization of~\eqref{eq:oneDimensionalExample} into $N$ degrees of freedom with $N+1$ uniform first-order Lagrange finite elements results in the system for the interior DOFs,
\begin{equation}\label{eq:oneDimensionalExampleDiscrete}
\mathbf{K} \displacementDiscreteFEM = \mathbf{b}
\end{equation}
where
$$
\stiffnessMatrix = \frac{AE}{\Delta x}
\begin{bmatrix}
2 & -1 &  \\
-1 & 2 & -1 &  \\
  & \ddots & \ddots  & \ddots &  \\
&  &  -1 & 2 & -1 \\
&  &  &  -1 & 2 \\
\end{bmatrix},
\qquad
\mathbf{b} = \Delta x 
\begin{bmatrix} 
 b \\  b  \\ \vdots \\ b \\ b
\end{bmatrix}
$$
with $\Delta x = \frac{L}{N+1}$. We note that the system is SPD, i.e., $\displacementDiscreteFEM^T \stiffnessMatrix \displacementDiscreteFEM > 0$, and as such is invertible such that the system~\eqref{eq:oneDimensionalExampleDiscrete} has a unique solution.

We now consider the domain decomposition formulation such that the first two and last two nodes are the ``coarse-scale" nodes and the remaining inner nodes are the ``fine-scale" nodes, i.e., $\displacementCoarseDiscreteFEM = \left[{\displacementDiscreteScalarFEM}_1,{\displacementDiscreteScalarFEM}_2,{\displacementDiscreteScalarFEM}_{N-1},{\displacementDiscreteScalarFEM}_{N} \right]^T \in \RR{4}$ and  $\displacementFineDiscreteFEM = \left[{\displacementDiscreteScalarFEM}_3,\cdots,{\displacementDiscreteScalarFEM}_{N-2} \right]^T \in \RR{N-4}$ with boundary degrees of freedom $\displacementCoarseDiscreteBoundaryFEM = \left[ {\displacementDiscreteScalarFEM}_2,{\displacementDiscreteScalarFEM}_{N-1} \right]^T \in \RR{2}$. The resulting coarse-scale problem is
\begin{equation}\label{eq:exampleCoarseScale}
\stiffnessMatrixCoarse \displacementCoarseDiscreteFEM + \closureTerm = \overline{\mathbf{b}} 
\end{equation}
where
$$\stiffnessMatrixCoarse =  
 \begin{bmatrix}
2 & -1 & 0 & 0 \\
-1 & 1 & 0 & 0 \\
0 & 0 & 1 & -1 \\
0 & 0 & -1 & 2 \\
\end{bmatrix}, 
\qquad 
\overline{\mathbf{b}} = 
\begin{bmatrix}
 b \\ b \\ b \\  b
\end{bmatrix},
$$
are the coarse-scale stiffness matrices and forcing vectors and
\begin{equation*}
 \closureTerm = 
\frac{AE}{\Delta x}
\begin{bmatrix}
0 & 0 & 0 & 0\\
1 & -1 & 0 & 0\\
0 & 0 & -1 &1\\
0 & 0 & 0 & 0\\
\end{bmatrix}
\begin{bmatrix}
{\displacementCoarseDiscreteScalarFEM}_2  \\ {\displacementFineDiscreteScalarFEM}_1 \\ {\displacementFineDiscreteScalarFEM}_{N-4} \\ {\displacementCoarseDiscreteScalarFEM}_{3} 
\end{bmatrix}
\end{equation*}
is the contribution of the interior degrees of freedom on the exterior degrees of freedom at the interface. We emphasize that this term depends on the (unknown) interior degrees of freedom. We further emphasize that the coarse-scale stiffness matrix is SPD such that $\displacementCoarseDiscreteFEM^T \stiffnessMatrixCoarse \displacementCoarseDiscreteFEM > 0$.

Replacing $\closureTerm$ from Eq.~\eqref{eq:exampleCoarseScale} with one of the machine-learned models $\mlModel$ outlined in Section~\ref{sec:ml-framework} results in the coarse-scale problem
 \begin{equation}\label{eq:exampleCoarseScaleb}
\stiffnessMatrixCoarse \displacementCoarseDiscreteFEM + \mlModelVec(\displacementCoarseDiscreteBoundaryFEM) = \overline{\mathbf{\forcingScalar}}. 
\end{equation}
We now provide a brief analysis of the models outlined in Section~\ref{sec:ml-framework}. 

\begin{itemize}
\item \textbf{Linear least-squares.} For this model the reduced system becomes
 \begin{equation}\label{eq:exampleCoarseScaleLls}
\left[ \stiffnessMatrixCoarse  + \stiffnessMatrix_{\mathsf{ML}} \right] \displacementCoarseDiscreteFEM  = \overline{\mathbf{\forcingScalar}} - \forceOffset 
\end{equation}
where
$$\stiffnessMatrix_{\mathsf{ML}} = 
\begin{bmatrix}
0 & 0 & 0 & 0 \\ 
0 & \left[ \basisForces \reducedLinearRegressionMatrix \basisDisplacements^T \right]_{11} & \left[ \basisForces \reducedLinearRegressionMatrix \basisDisplacements^T \right]_{12}   & 0 \\
0 & \left[ \basisForces \reducedLinearRegressionMatrix \basisDisplacements^T \right]_{21} & \left[ \basisForces \reducedLinearRegressionMatrix \basisDisplacements^T \right]_{22}   & 0 \\
0 & 0 & 0 & 0\\
\end{bmatrix}.
$$
As we place no constraints on the construction of $\reducedLinearRegressionMatrix$, we cannot guarantee that $\stiffnessMatrix_{\mathsf{ML}}$ is positive (semi-)definite, nor can we guarantee a unique solution to the system~\eqref{eq:exampleCoarseScaleLls} without knowledge of the operator  $\basisForces \reducedLinearRegressionMatrix \basisDisplacements^T$, which in general is not full rank. As a canonical example, we can consider the case where 
$$
\basisForces \reducedLinearRegressionMatrix \basisDisplacements^T = 
\begin{bmatrix} -0.5& 0 \\ 0 & -0.5 \end{bmatrix},
$$
which is an admissible solution to the optimization problem~\eqref{eq:reducedLeastSquaresMatrix}. In this case the coarse-scale equation stiffness matrix becomes
$$\stiffnessMatrixCoarse =  
 \begin{bmatrix}
2 & -1 & 0 & 0 \\
-1 & 0.5 & 0 & 0 \\
0 & 0 & 0.5 & -1 \\
0 & 0 & -1 & 2 \\
\end{bmatrix}
$$
which is singular. 
We further observe that, in general, the ML contribution will \textit{not} be symmetric. Thus, even if invertible, the system~\eqref{eq:exampleCoarseScaleLls} cannot in general be solved with the popular conjugate gradient method. 

\item \textbf{Neural network.} The neural network model results in the coarse-scale equation
 \begin{equation}\label{eq:exampleCoarseScaleNn}
 \stiffnessMatrixCoarse  \displacementCoarseDiscreteFEM  +  \basisForces \reducedNeuralNetwork \left( \basisDisplacements^T \displacementCoarseDiscreteFEM^{\domainBoundaryOneTwo} ; \weights \right) = \overline{\mathbf{\forcingScalar}} - \forceOffset. 
\end{equation}
Similar to the linear least-squares model, we cannot guarantee that the Newton iteration associated with~\eqref{eq:exampleCoarseScaleNn} will involve full-rank matrices. In general the matrices will not be SPD and solution via the conjugate gradient method is inappropriate. 

\item \textbf{Symmetric positive semi-definite linear-least-squares} The SPSD-LLS model results in the reduced system
 \begin{equation}\label{eq:exampleCoarseScaleSpdLls}
\left[ \stiffnessMatrixCoarse  + \stiffnessMatrix_{\mathsf{ML}} \right] \displacementCoarseDiscreteFEM  = \overline{\mathbf{\forcingScalar}} - \forceOffset 
\end{equation}
where
$$\stiffnessMatrix_{\mathsf{ML}} = 
\begin{bmatrix}
0 & 0 & 0 & 0 \\ 
0 & \left[ \basisCombined \reducedLowerLinearRegressionMatrix \reducedLowerLinearRegressionMatrix^T \basisCombined^T \right]_{11}   & \left[ \basisCombined \reducedLowerLinearRegressionMatrix \reducedLowerLinearRegressionMatrix^T \basisCombined^T \right]_{12}   & 0 \\
0 &\left[ \basisCombined \reducedLowerLinearRegressionMatrix \reducedLowerLinearRegressionMatrix^T \basisCombined^T \right]_{21} & \left[ \basisCombined \reducedLowerLinearRegressionMatrix \reducedLowerLinearRegressionMatrix^T \basisCombined^T \right]_{22}   & 0 \\
0 & 0 & 0 & 0\\
\end{bmatrix}.
$$
By construction, $\displacementCoarseDiscreteFEM^T \stiffnessMatrix_{\mathsf{ML}} \displacementCoarseDiscreteFEM \ge 0$ and, as a result, $\displacementCoarseDiscreteFEM^T\left[\stiffnessMatrixCoarse +  \stiffnessMatrix_{\mathsf{ML}}\right] \displacementCoarseDiscreteFEM > 0$. Thus, the coarse-scale system is SPD. By definition of SPD matrices, the system is invertible with a unique solution and may be solved with the conjugate gradient method.

\item \textbf{Symmetric positive semi-definite neural network}. The SPSD-NN model results in the reduced system
 \begin{equation}\label{eq:exampleCoarseScaleSpdLls}
\left[ \stiffnessMatrixCoarse  + \stiffnessMatrix_{\mathsf{ML}}(\displacementCoarseDiscreteFEM) \right] \displacementCoarseDiscreteFEM  = \overline{\mathbf{\forcingScalar}} - \forceOffset 
\end{equation}
where
$$\stiffnessMatrix_{\mathsf{ML}} = 
\begin{bmatrix}
0 & 0 & 0 & 0 \\ 
0 &\left[ \basisCombined \reducedNeuralNetworkStiffness  \reducedNeuralNetworkStiffness^T  \basisCombined^T  \right]_{11} & 
\left[ \basisCombined \reducedNeuralNetworkStiffness  \reducedNeuralNetworkStiffness^T  \basisCombined^T  \right]_{12}  & 0 \\
0 & \left[ \basisCombined \reducedNeuralNetworkStiffness  \reducedNeuralNetworkStiffness^T  \basisCombined^T  \right]_{21} & \left[ \basisCombined \reducedNeuralNetworkStiffness  \reducedNeuralNetworkStiffness^T  \basisCombined^T  \right]_{22}  & 0 \\
0 & 0 & 0 & 0\\
\end{bmatrix}.
$$
Again, by construction, $\displacementCoarseDiscreteFEM^T \stiffnessMatrix_{\mathsf{ML}} \displacementCoarseDiscreteFEM \ge 0$ and, as a result, the coarse-scale system is symmetric positive definite for any $\displacementCoarseDiscreteFEM$. The system is invertible with a unique solution and may be solved with the conjugate gradient method.
\end{itemize}

In summary, this section demonstrated that, in the simplified setting of linear elasticity in one-dimension, the SPSD stiffness-based models result in coarse-scale systems that are SPD. These systems have invertible system matrices with unique solutions. We further demonstrated that we cannot guarantee the systems emerging from the direct displacement-to-force mappings have unique solutions. 
\section{Numerical examples}\label{sec:experiment}
We now investigate the performance of the proposed models on several exemplars, starting from a canonical cube undergoing finite deformations to a more complex fastener-bushing geometry undergoing finite deformations with contact.
  
\subsection{Implementation in Sierra Solid Mechanics}
The machine-learned models outlined in Section~\ref{sec:ml-framework} are implemented in the \textsf{Sierra/SolidMechanics} code base,
  which is a part of the \textsf{Sierra} simulation code suite developed at Sandia National Laboratories.
The ML implementation is built on top of the \textsf{pocket-tensor library} \cite{pocket-tensor},
  which provides a C++ interface for \textsf{TensorFlow} and \textsf{Keras} models.
For the direct displacement-to-force models,
  the stiffness matrix associated with the ML element is computed via finite differences,
  while for the displacement-to-stiffness models the stiffness matrix is directly output from the neural network.
For all cases, the nonlinear system of equations resulting from the FEM and FEM--ML models is solved via the nonlinear conjugate gradient method with a full tangent preconditioner. We employ hyperelastic material models for all examples and solve the governing equations with finite deformations under a quasi-static approximation. Contact is handled using the Augmented Lagrange algorithm. 
We refer the interested reader to the \textsf{Sierra/SolidMechanics} theory manual for more details~\cite{sierra}.
We note that, while it would be interesting to assess the performance of more generic solvers for the non-SPSD models (e.g., generalized minimum residual),
  these are not available in \textsf{Sierra/SolidMechanics} given the superior performance of conjugate gradient for SM problems.

\subsection{Training of ML models}\label{sec:ml_training}
All machine-learned models are trained in \textsf{Python}. For LLS, we directly solve the least-squares problem using \textsf{scipy.optimize.lsq\_linear}. Given the analytic solution to the optimization problem we do not hold out any data in training. For SPSD-LLS, we solve the optimization problem with the splitting conic solver as implemented in the software \textsf{cvxpy}~\cite{diamond2016cvxpy,agrawal2018rewriting}. We again do not hold out any data for training SPSD-LLS. 

Both the NN and SPSD-NN models are implemented and trained in \textsf{PyTorch}. For training the networks we employ random shuffling and a standard 80/20 training/validation split of the data. The neural networks are trained with the ADAM optimizer for 35000 epochs. We employ an early stopping criteria if the 200-epoch running mean of the validation loss has not decreased. Unless otherwise noted, we employ a batch size of $500$ with a learning rate schedule of $\textsf{lr} = \{1 \times 10^{-3}, 2\times 10^{-4}, 1 \times 10^{-4} , 5 \times 10^{-5}, 2 \times 10^{-5} \}$, where the learning rate is lowered after $500, 1000, 2000, 5000$, and $15000$ epochs. For network architecture, we employ fully connected neural networks with three hidden layers. The number of nodes in each layer is set to be equivalent to the size of the network input vector (e.g., the reduced basis dimension, $K$). The loss function for all models is the $\ell^2$ loss between the predicted and truth forces. Lastly, we employ ReLU for our nonlinear activation function. PyTorch code for both the standard neural network architecture and the SPSD neural network architecture are provided in~\ref{sec:pytorch_code}.

The training stage for our examples involves solving a full FEM that directly resolves the fastener. The coarse and fine scales are then taken to be a subset of this full FEM discretization and the training data are extracted accordingly. In theory, one could employ non-conforming decompositions (e.g., perform the training on one discretization and learn an ML model that is applicable for other meshes via, e.g., interpolation). For simplicity, we do not consider this here. In practice, our training process consists of two phases. First, we solve the full FEM that directly resolves the coarse and fine scales. After this, we extract the fine-scale domain from the model and execute a post-processing step to compute the internal forces. 
\subsection{Multi-element deformed cube}
\begin{figure}
\begin{subfigure}[t]{0.45\textwidth}
\includegraphics[trim={6cm 2cm 6cm 2cm},clip,width=1.0\linewidth]{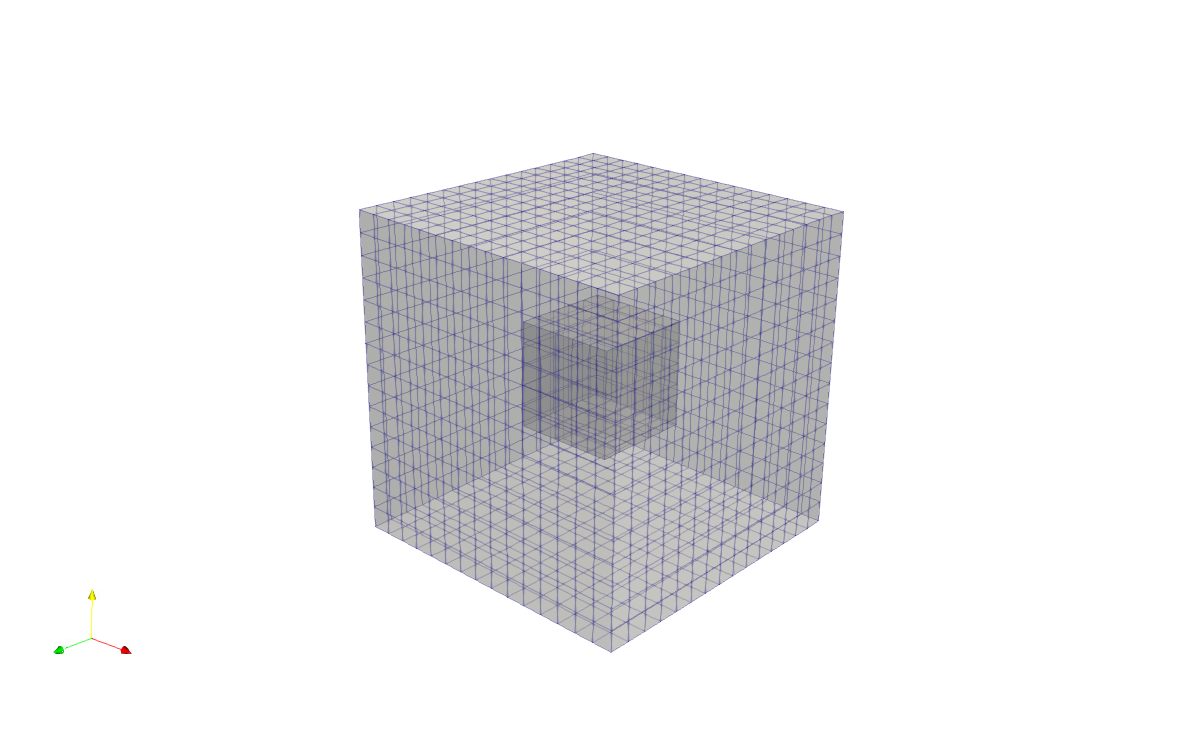}
\end{subfigure}
\begin{subfigure}[t]{0.45\textwidth}
\includegraphics[trim={5cm 3cm 7cm 2cm},clip,width=1.0\linewidth]{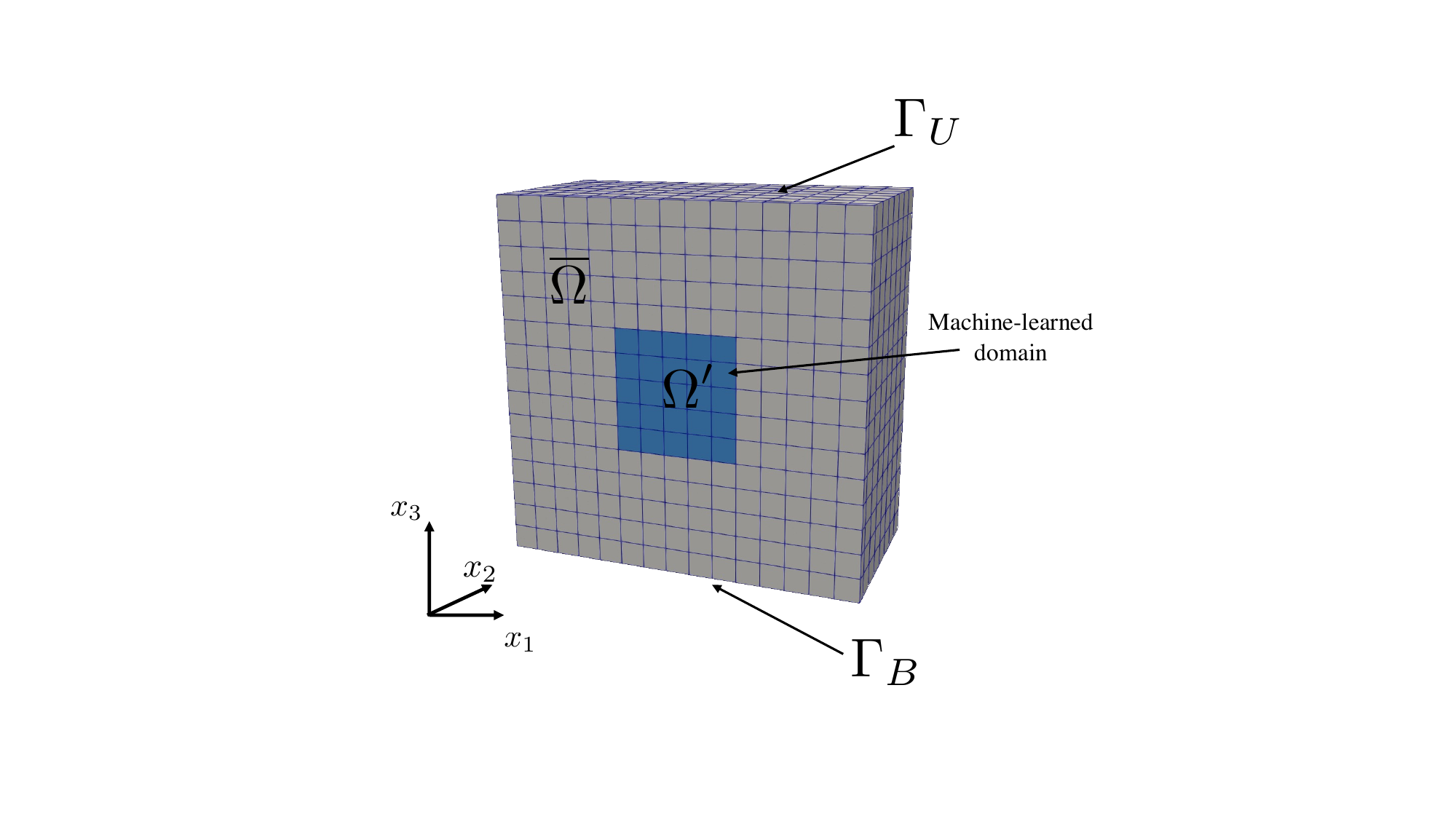}
\end{subfigure}
\caption{\cubeCaption\ Depiction of the setup of the cube problem. The ML model replaces the interior cube, which comprises a $5\times5\times 5$ element block of the full FEM model.}
\label{fig:cubeSetup}
\end{figure}

We begin by considering quasi-static deformation of a cube on the domain $\domain \equiv [-1.5 \mathsf{ in},1.5 \mathsf{ in}]^3$. The ``inner" fine-scale domain is defined as $\domainTwo  = (-0.5 \mathsf{ in },0.5 \mathsf{ in})^3$ while the outer domain is defined as $\domainOne = \domain - \domainTwo$. For discretization, the domain $\domain$ is discretized into  $15 \times 15 \times 15$ uniformly sized cubes in each direction resulting in $\nCoarseNodesDomainBoundaryOneTwo = 152$ interface nodes. We employ standard bilinear Lagrange finite elements. Figure~\ref{fig:cubeSetup} depicts a setup of the problem. 
We employ homogeneous Dirichlet boundary conditions on the bottom boundary and Dirichlet boundary conditions on the top boundary, i.e.,
\begin{equation}
\begin{split}
\displacementFEM &= \mathbf{0} \text{ on } \Gamma_B ,\\
\displacementFEM &= \left[ u_1, u_2,u_3 \right] \text{ on } \Gamma_U,
\end{split}
\end{equation}
where  $\Gamma_{B} = [-1.5 \mathsf{ in},1.5 \mathsf{ in}]  \times [-1.5  \mathsf{ in},1.5  \mathsf{ in}] \times -1.5  \mathsf{ in}$ is the lower boundary and $\Gamma_{U} = [-1.5  \mathsf{ in},1.5  \mathsf{ in}]\times [-1.5  \mathsf{ in},1.5  \mathsf{ in} ]  \times 1.5  \mathsf{ in}$ is the upper boundary. We employ a hypoelastic material model across the entire domain. For $\x \in \domainOne$ the material model is characterized by a Young's modulus of 28.5e6 \textsf{psi} with a Poisson's ratio of $\nu = 0.3$. For $\x \in \domainTwo$ the material model is characterized by a Young's modulus of $29 \times 10^6$ \textsf{psi} with a Poisson's ratio of $\nu=0.3$. As a quantity of interest (QoI), we consider integrated reaction forces on the bottom of the cube, $\Gamma_B$. 

\begin{figure}
\begin{center}
\begin{subfigure}[t]{0.45\textwidth}
\includegraphics[trim={0cm 0cm 0cm 0cm},clip,width=1.0\linewidth]{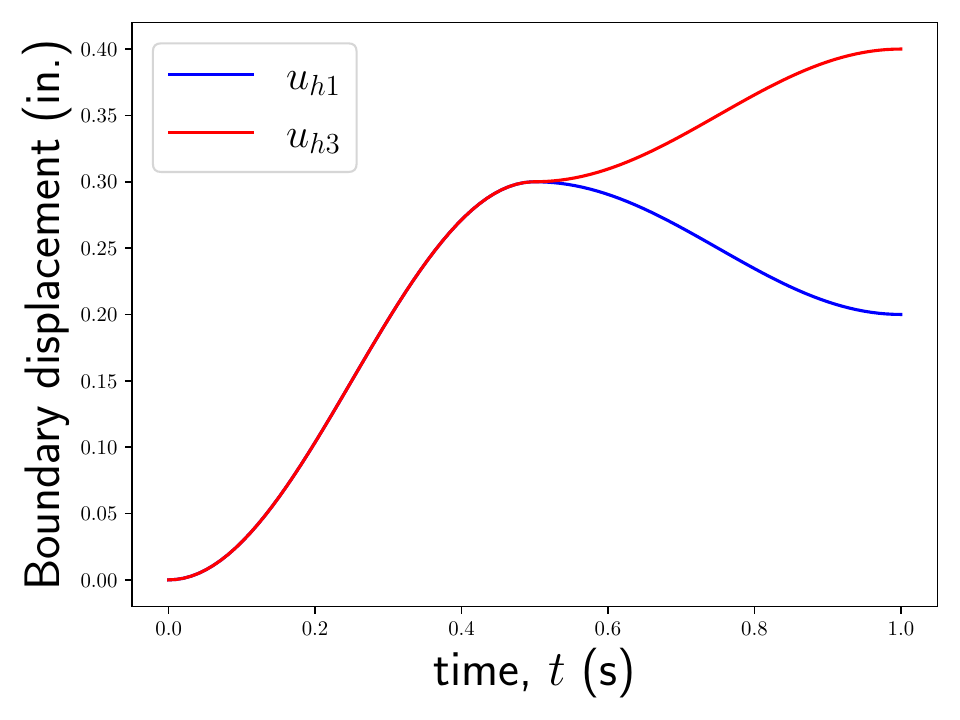}
\end{subfigure}
\caption{\cubeCaption\ Depiction of loading profile for $\scaleFactorXOne = 0.3$, $\scaleFactorXTwo=-0.1$, $\scaleFactorZOne=0.3$, $\scaleFactorZTwo=0.1$.}
\label{fig:cubeLoading}
\end{center}
\end{figure}

For training we solve the cube problem for 16 different non-proportional quasi-static loading trajectories for $t \in [0\mathsf{s},1\mathsf{s}]$, where $t$ is pseudo-time and a single loading trajectory comprises 100 time steps with the top boundary displaced by a fixed profile. This loading profile is prescribed (in inches) for $\x \in \Gamma_U$ as 
\begin{equation}
\begin{split}
&{\displacementFEM}_{1}(t,\x) = \begin{cases}
 \cosRamp{\scaleFactorXOne}{0}{\timeTermXOne}   \qquad & 0 \le t \le \timeTermXOne  \\
 \scaleFactorXOne + \cosRamp{\scaleFactorXTwo}{\timeTermXOne}{T} \qquad & \timeTermXOne < t < T \\
\end{cases} \\
&{\displacementFEM}_{3}(t,\x) = \begin{cases}
 \cosRamp{\scaleFactorZOne}{0}{\timeTermZOne}   \qquad & 0 \le t \le \timeTermZOne\\
 \scaleFactorZOne + \cosRamp{\scaleFactorZTwo}{\timeTermZOne}{T} \qquad & \timeTermZOne < t < T\\
\end{cases}
\end{split}
\end{equation}
where $\timeTermXOne = \timeTermZOne = 0.5\mathsf{s}$. Figure~\ref{fig:cubeLoading} shows an example loading profile for $\scaleFactorXOne = 0.3$, $\scaleFactorXTwo=-0.1$, $\scaleFactorZOne=0.3$, $\scaleFactorZTwo=0.1$.

As a training set we examine a hyper-cube of the parameters $\scaleFactorXOne = (-0.3,0.3)$, $\scaleFactorXTwo = (-0.1,0.1)$, $\scaleFactorZOne = (-0.3,0.3)$, and $\scaleFactorZTwo = (-0.1,0.1)$. In all, this results in 1600 total quasi-static solves (16 parameter configurations with 100 time steps per configuration). For testing we examine four loading configurations as described in Table~\ref{tab:cubeTestingParameters}.  We note that this setup comprises training runs that deform the body by up to 10\% and testing configurations that deform the body up to 16\%, both of which are beyond the small-strain region and result in non-linearities.  

\begin{table}
\begin{tabular}{ |c|c|c| c| c|} 
 \hline
 Testing run number & $\scaleFactorXOne$ & $\scaleFactorXTwo$ & $\scaleFactorZOne$ & $\scaleFactorZTwo$ \\
 \hline \hline
 1 & 0.05 & -0.01 & -0.02 & 0.01\\  \hline
 2 & 0.15 & -0.10 & -0.20 & 0.15 \\  \hline
 3 & 0.40 & -0.10 & -0.30 & 0.10 \\ \hline
 4 & -0.30 & -0.30 & 0.50 & -0.05 \\ \hline
\end{tabular}
\caption{\cubeCaption\ Parameter values for testing configurations.}
\label{tab:cubeTestingParameters}
\end{table}

\subsubsection{Offline training results}
\begin{figure}
\begin{center}
\begin{subfigure}[t]{0.49\textwidth}
\includegraphics[trim={0cm 0cm 0cm 0cm},clip,width=1.0\linewidth]{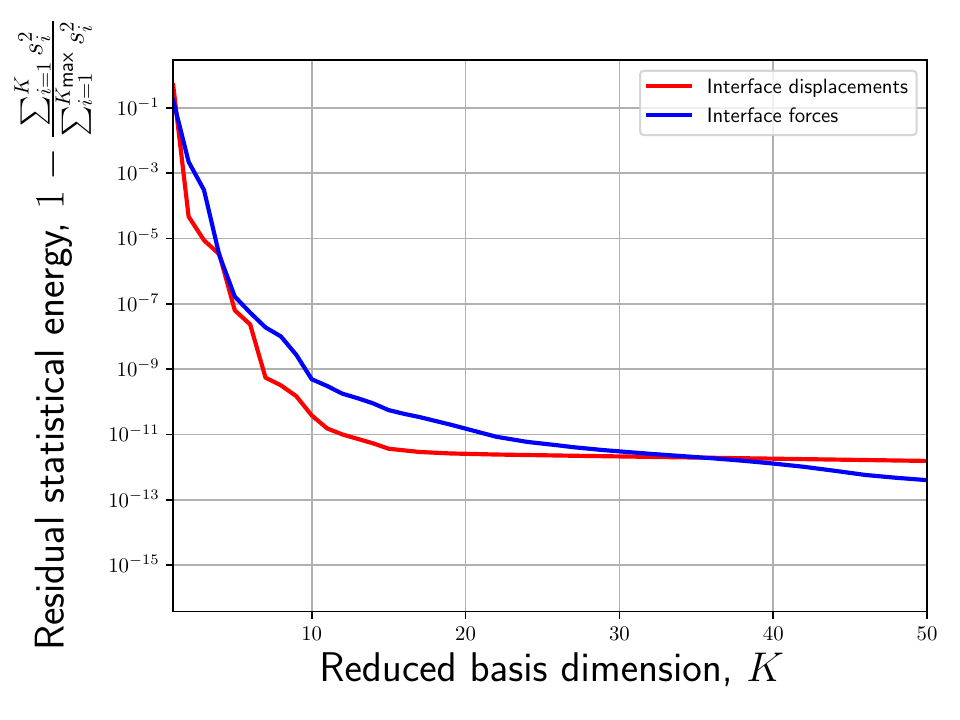}
\end{subfigure}
\caption{\cubeCaption\ Residual statistical energy of the reduced basis approximation for the interface displacement and force fields.}
\label{fig:cubeTrainingBasisConvergence}
\end{center}
\end{figure}

We first assess the reduced basis approximations as described in the optimization problems~\eqref{eq:basis_force} and~\eqref{eq:basis_displacement}. As there is no body forcing term, we take $\forceOffset = \mathbf{0}$. Figure~\ref{fig:cubeTrainingBasisConvergence} depicts the residual statistical energy (i.e., relative energy contained in the truncated singular values) associated with the reduced basis approximation\footnote{Note that the residual statistical energy bounds the error of the reduced basis approximation.}. We observe that both the interface force field and interface displacement field are well-represented with relatively few basis vectors. With just 10 basis vectors the residual statistical energy of both the interface displacement and force fields is less than $10^{-9}$. This quick decay demonstrates that the interface can be characterized with relatively few degrees of freedom, and justifies the ansatz of performing the learning process in the reduced space. 

\begin{figure}
\begin{center}
\begin{subfigure}[t]{0.49\textwidth}
\includegraphics[trim={0cm 0cm 0cm 0cm},clip,width=1.0\linewidth]{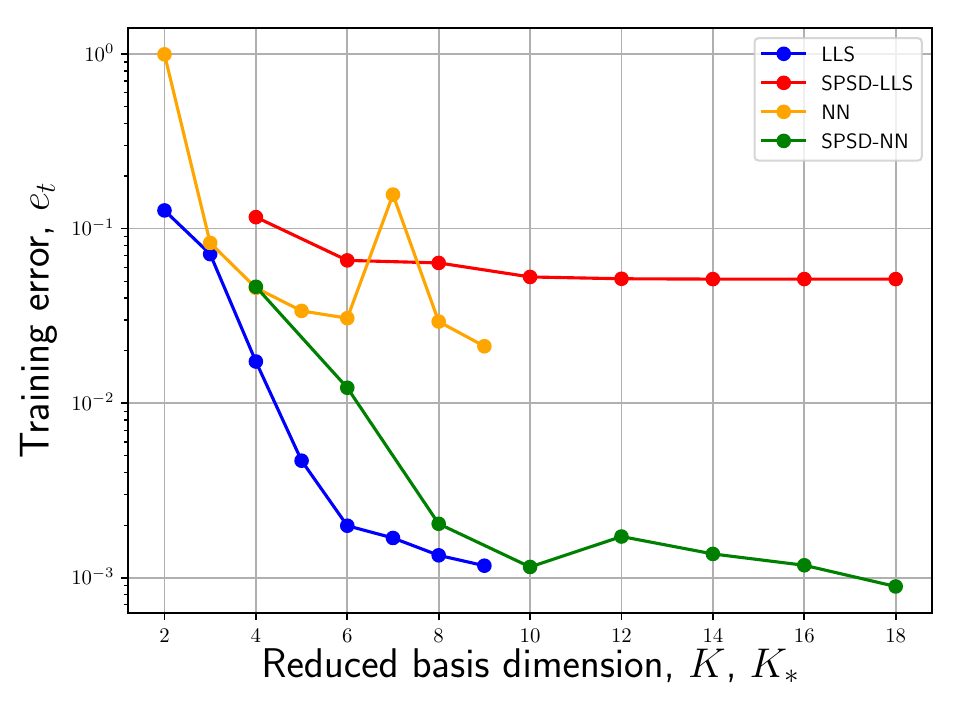}
\end{subfigure}
\begin{subfigure}[t]{0.49\textwidth}
\includegraphics[trim={0cm 0cm 0cm 0cm},clip,width=1.0\linewidth]{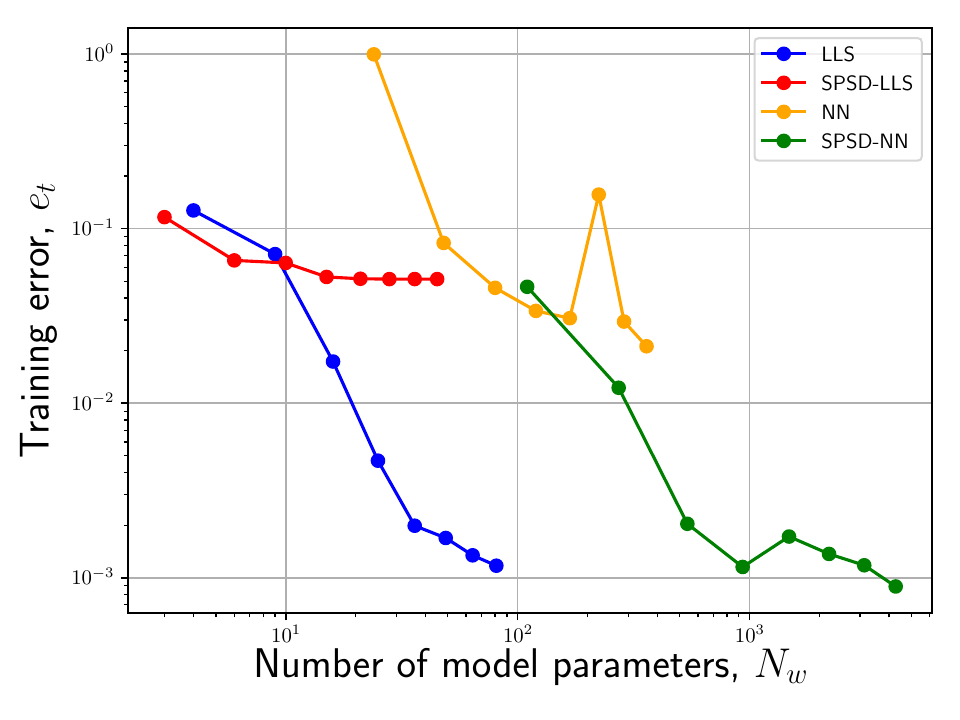}
\end{subfigure}
\caption{\cubeCaption\ Offline training results of machine learned models. On the left we show the training error as a function of the reduced basis dimension $K$ while on the right we show results as a function of model parameters, $N_w$. Note that the SPSD models employ the combined basis of dimension $K^* = 2K$.}
\label{fig:cubeTrainingResults}
\end{center}
\end{figure}

Next, Figure~\ref{fig:cubeTrainingResults} shows the training error for the various ML models considered as a function of the reduced basis dimension, $K$, and the number of model parameters, $N_w$. We observe that LLS yields significantly more accurate results than SPSD-LLS. SPSD-LLS saturates in accuracy around $\reducedDimensionCombined = 10$ while LLS continues to converge, reaching a training error of $10^{-3}$. Surprisingly, we next observe that LLS results in more accurate models than NN. In theory, NN should be at least as accurate as LLS. We believe that, for this example, LLS performs better than NN because the data are almost linear (recall we employ a hypoelastic material model). Since LLS has an analytic solution, we are more easily able to optimize the model. We additionally note that we observe a ``kink" in the convergence of NN at $K=7$. We believe that this kink is a result of the stochastic training of the NNs, and we could likely eliminate it by training an ensemble of models and employing cross-validation to select the best model (presently we train one model with early stopping, as described in Section~\ref{sec:ml_training}). Lastly, we observe that SPSD-NN significantly outperforms both NN and SPSD-LLS. We expect that SPSD-NN greatly outperforms NN due to the structure embedded in the model form. In particular, we are learning an SPSD stiffness matrix which is similar to the structure one would expect from theory.

\subsubsection{Coupled FEM-ML results}
\begin{figure}
\begin{center}
\begin{subfigure}[t]{0.48\textwidth}
\includegraphics[trim={0cm 0cm 0cm 0cm},clip,width=1.0\linewidth]{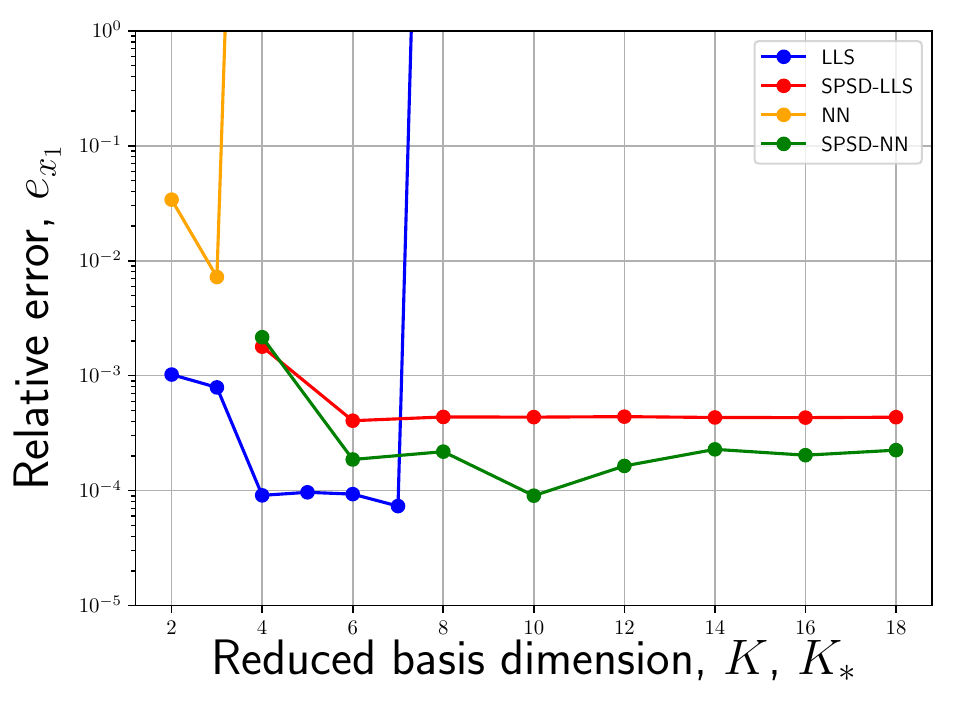}
\end{subfigure}
\begin{subfigure}[t]{0.48\textwidth}
\includegraphics[trim={0cm 0cm 0cm 0cm},clip,width=1.0\linewidth]{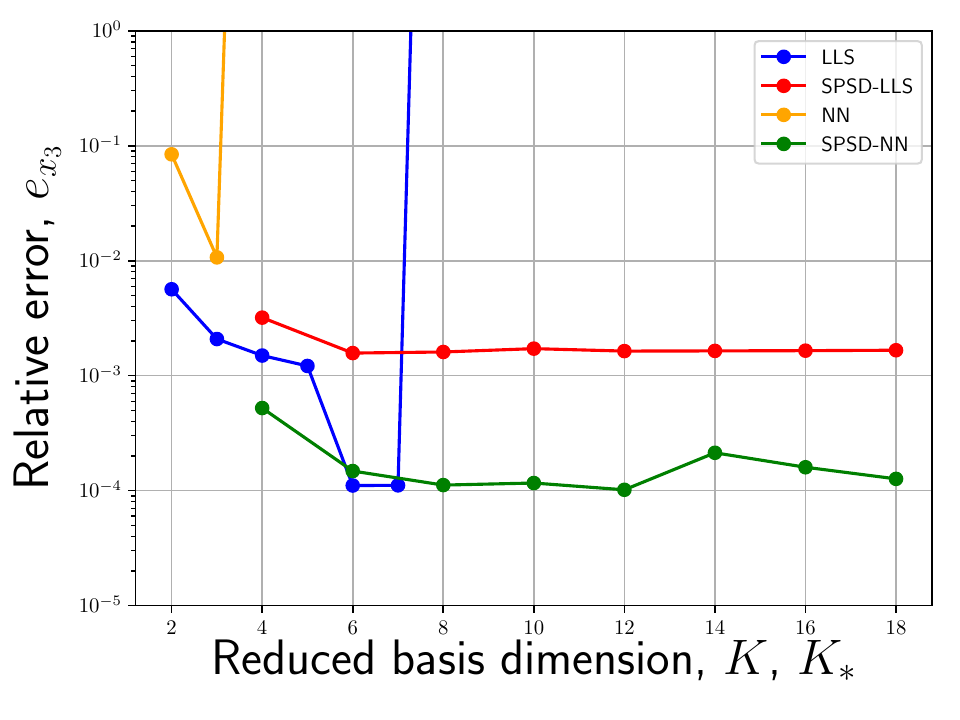}
\end{subfigure}
\caption{\cubeCaption\ Relative errors for coupled FEM-ML models as a function of reduced basis dimension. Note that the SPSD models employ the combined basis of dimension $K^* = 2K$.}
\label{fig:cubeApostConvergence}
\end{center}
\end{figure}

We now examine \textit{a posteriori} results where the ML model is coupled to the FEM solver. Figure~\ref{fig:cubeApostConvergence} shows QoI errors for the various models as a function of the reduced dimension. LLS results in accurate models for small reduced basis dimensions but quickly goes unstable as the basis dimension grows past $\reducedDimension = 7$. We emphasize that this instability is not observed in the offline training results and is likely a result of the interplay between the LLS-POD model and the solver. Despite being the least accurate model in training, SPSD-LLS results in stable and accurate models with relative errors for both QoIs of less than 1\% for all reduced basis dimensions. NN results in inaccurate models that quickly go unstable for reduced basis dimensions $\reducedDimension > 3$. SPSD-NN is the most accurate and robust model. It is stable for all reduced basis dimensions and results in relative errors of less than  $0.1\%$ for $\reducedDimensionCombined \ge 6$.

\begin{figure}
\begin{center}
\begin{subfigure}[t]{0.95\textwidth}
\includegraphics[trim={0cm 0cm 0cm 0cm},clip,width=1.0\linewidth]{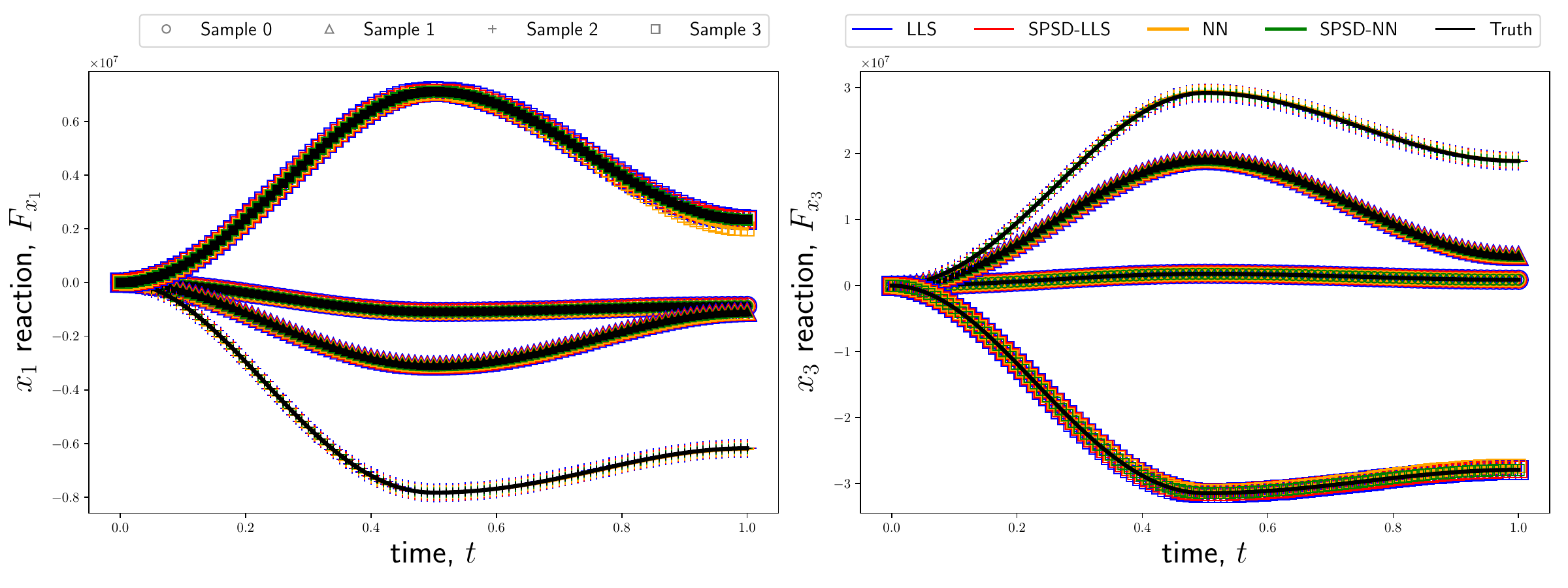}
\end{subfigure}
\caption{\cubeCaption\ Results of ML--FEM and FEM-only model on all testing configurations. Results are shown for the optimal configuration for each model, which corresponds to $K=7$, $K^*=14$, $K=3$, and $K^* =10$ for LLS, SPSD-LLS, NN, and SPSD-NN, respectively.}
\label{fig:cubeProfile}
\end{center}
\end{figure}

\begin{figure}
\begin{center}
\begin{subfigure}[t]{1.0\textwidth}
\includegraphics[trim={2cm 8cm 4cm 10cm},clip,width=1.0\linewidth]{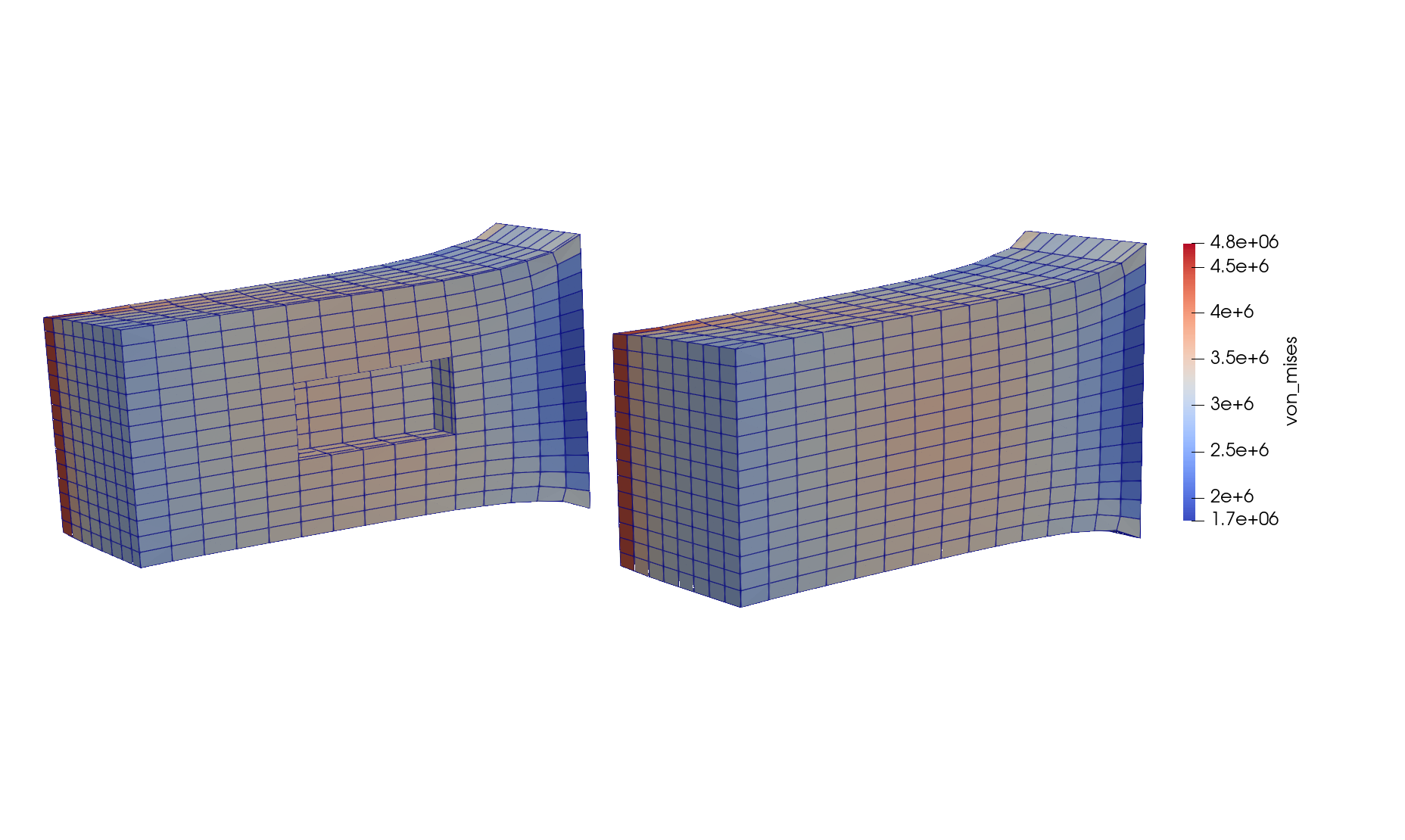}
\end{subfigure}
\caption{\cubeCaption\ Results of FEM-ML coupled simulation with the SPSD-NN model (left) and FEM-only simulation (right) for prediction of the Von-Mises stress. Note that displacements are amplified for visualization.}
\label{fig:cubeVonMises}
\end{center}
\end{figure}

Next, Figure~\ref{fig:cubeProfile} shows QoI results for all testing configurations. We plot results for the optimal configurations (i.e., the configuration leading to the lowest error) for each model. The optimal configurations are $K=7$, $K^*=14$, $K=3$, and $K^* =10$ for LLS, SPSD-LLS, NN, and SPSD-NN, respectively. We observe that all predictions lie on top of the truth data with the exception of the NN model, which deviates slightly from the truth near $t=1.0$ for testing configuration number 3. Overall, we observe that all models are accurate when they converge. Lastly, Figure~\ref{fig:cubeVonMises} shows the von Mises stress predicted by the SPSD-NN model (left) and full model (right) for training configuration \#4. With the contours shown in the figure, we observe no noticeable difference between the FEM--ML and FEM-only solutions. 

The results of this elementary example demonstrate the capability of the ML models to replace a subdomain, even when the behavior of the material is in a nonlinear large strain regime. We additionally observed that, while LLS and NN models were accurate in the offline phase, they rarely led to converged solutions when coupled with the FEM solver. The SPSD-LLS and SPSD-NN models, both of which enforced an SPSD stiffness matrix, were more robust when coupled to the FEM model. 


\subsection{Fastener undergoing loading with contact}
\begin{figure}
\begin{center}
\begin{subfigure}[t]{1.0\textwidth}
\includegraphics[trim={4cm 1cm 4cm 2cm},clip,width=1.0\linewidth]{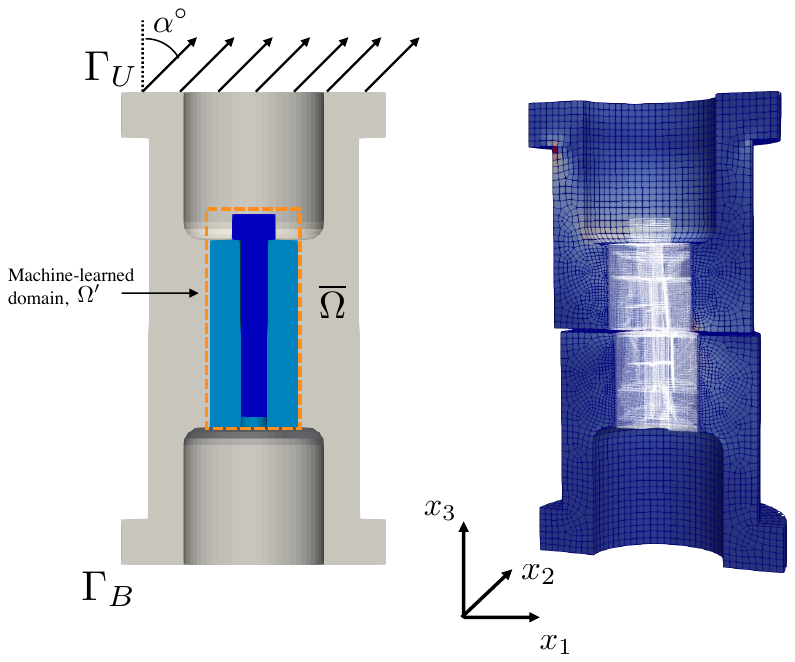}
\end{subfigure}
\caption{\contactCaption\ Problem schematic for fastener example. On the left we show the decomposition between the resolved FEM domain and the ML domain along with the loading configuration, while on the right we show a sample training solution at $\alpha \approx 65^{\circ}$. We show contours of max Von Mises stress with the fastener subdomain colored in white. Displacements are amplified for visualization.}
\label{fig:fastenerSetup}
\end{center}
\end{figure}

We now consider a more complex example comprising a fastener-bushing geometry undergoing quasi-static radial loading. 
This geometry is derived from a test setup which can be found in~\cite{MoAvMi21}. The NAS1351-3-20P fastener is modeled as a ``plug" (fastener head and shank with no modeling of threads) with a 0.187 inch diameter and sits in the middle of the top bushing through hole, initially not in contact. The fastener is connected to the bottom bushing with a contiguous mesh. As the joint is loaded laterally, the 0.003 inch gap between the top bushing and fastener will close and the two volumes will be in contact. The problem is symmetric about the $x_3$-axis, and symmetry is enforced with a symmetric boundary condition. There is no preload in this case; preload will be considered in the next section. Figure~\ref{fig:fastenerSetup} depicts the problem. The full FEM mesh has $N=115,292$ degrees of freedom and only comprises half the domain. We remove both the fastener and the immediate domain around the fastener, which results in an interface with $\nCoarseNodesDomainBoundaryOneTwo = 3348$ nodes. We emphasize that contact occurs exclusively within the removed domain, and as a result our ML model must be able to accurately characterize contact. We employ homogeneous Dirichlet conditions on the exterior of the bottom bushing and Dirichlet boundary conditions on the top boundary. The quasi-static loading profile is given in inches for $\x$ on $\Gamma_U$, $t \in [0,1]$ by ${\displacementFEM}_{1}(t,\x) =  \beta t$, ${\displacementFEM}_3(t,\x) = \frac{1}{400} t$; this loading configuration corresponds to radial loading (i.e., a pull at a constant angle) whose angle is a function of the parameter $\beta$. 
We again employ a hypoelastic material model across the entire domain. For $\x \in \domainOne$ the material model is characterized by a Young's modulus of 28.5e6 \textsf{psi} with Poisson's ratio of $\nu = 0.3$. For $\x \in \domainTwo$ the material model is characterized by a Young's modulus of $29 \times 10^6$ \textsf{psi} with a Poisson's ratio of $\nu=0.3$. As QoIs, we consider the intergrated reaction forces on the exterior of the bottom bushing in the $x_1$ and $x_3$ directions. 
 
We consider 21 training configurations for $\beta = \{-0.005 + 0.0005 i\}_{i=0}^{20}$ and $t \in \{0.01i\}_{i=1}^{101}$. This dataset comprises quasi-static trajectories with 101 time steps per trajectory for pulls ranging from approximately $-65^{\circ}$ to $65^{\circ}$, which results in 2121 total quasi-static solves. 

\subsubsection{Offline training results}
\begin{figure}
\begin{center}
\begin{subfigure}[t]{0.49\textwidth}
\includegraphics[trim={0cm 0cm 0cm 0cm},clip,width=1.0\linewidth]{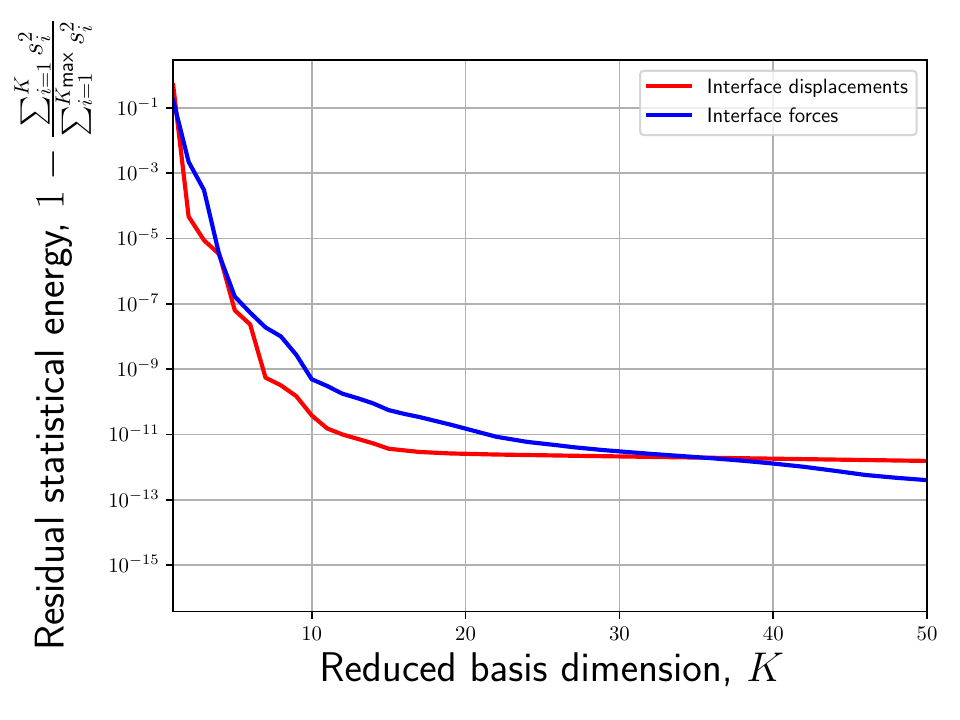}
\end{subfigure}
\caption{\contactCaption\ Residual statistical energy of the reduced basis approximation for the interface displacement and force fields.}
\label{fig:fastenerTrainingBasisConvergence}
\end{center}
\end{figure}

Figure~\ref{fig:fastenerTrainingBasisConvergence} demonstrates that the interface degrees of freedom are amenable to dimension reduction by displaying the residual statistical energy associated with the reduced basis approximation. As there is again no body forcing term, we take $\forceOffset = \mathbf{0}$.
As compared to the previous exemplar, the residual statistical energy in the force approximation decays slower. Overall, however, we again observe that both the interface force field and interface displacement field are well-represented with relatively few basis vectors. With just 10 basis vectors the residual statistical energy of the interface displacement field is less than $10^{-10}$ and the interface force field is less than $10^{-6}$. This quick decay again demonstrates that the interface can be characterized with relatively few dimensions, and justifies the ansatz of performing the learning process in the reduced space. 

\begin{figure}
\begin{center}
\begin{subfigure}[t]{0.49\textwidth}
\includegraphics[trim={0cm 0cm 0cm 0cm},clip,width=1.0\linewidth]{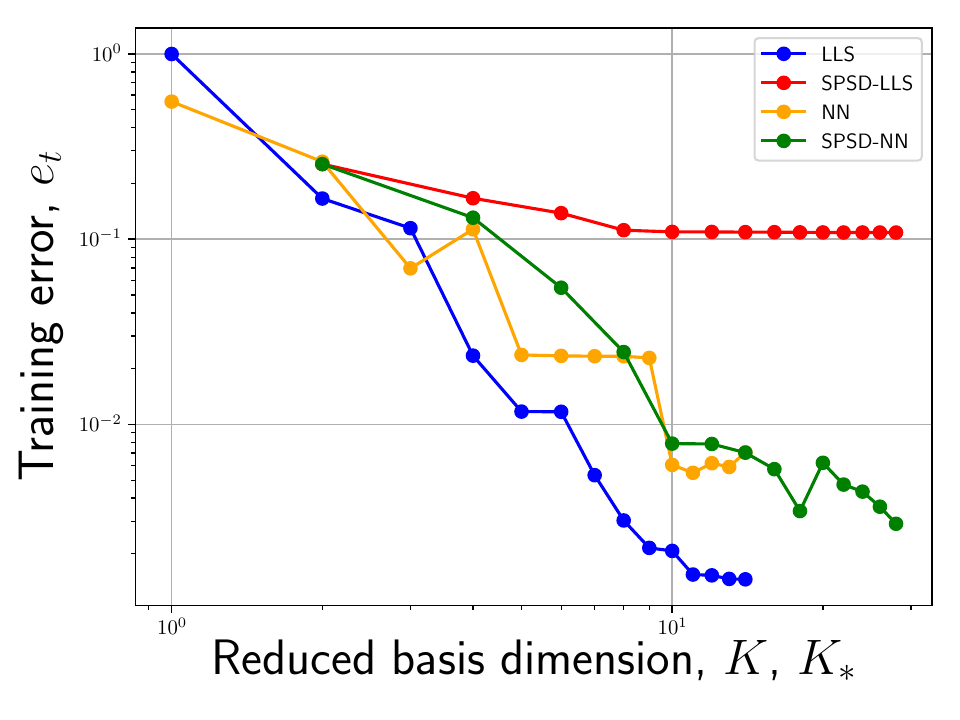}
\end{subfigure}
\begin{subfigure}[t]{0.49\textwidth}
\includegraphics[trim={0cm 0cm 0cm 0cm},clip,width=1.0\linewidth]{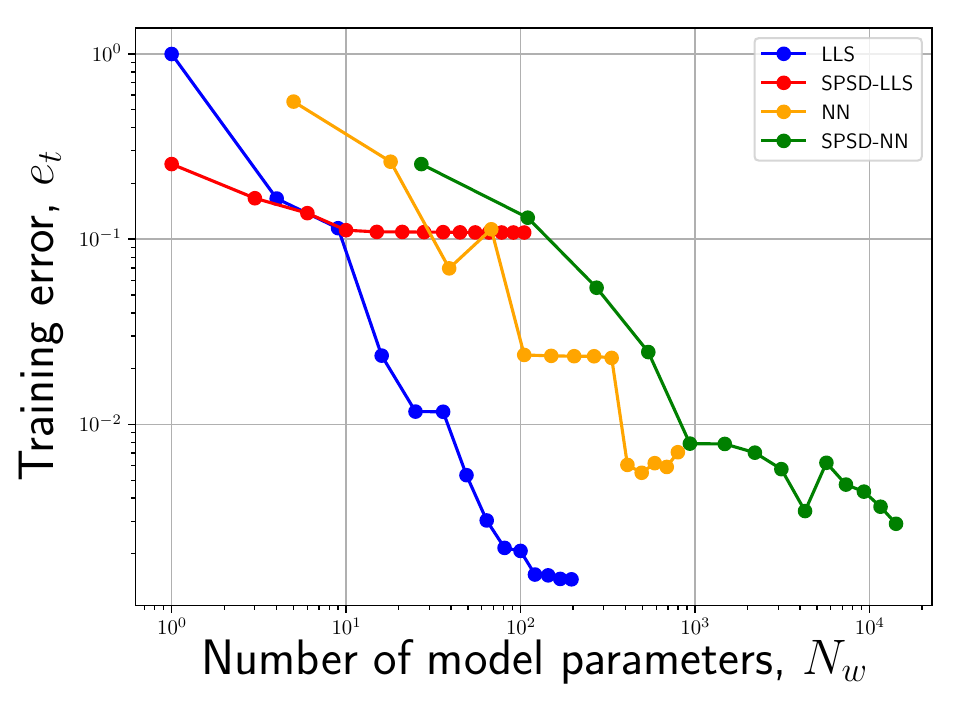}
\end{subfigure}
\caption{\contactCaption\ Offline training results of machine learned models trained on the fastener dataset. On the left we show the training error as a function of the reduced basis dimension $K$ while on the right we show results as a function of model parameters, $N_w$. Note that the SPSD models employ the combined basis of dimension $K^* = 2K$.}
\label{fig:fastenerTrainingResults}
\end{center}
\end{figure}

Figure~\ref{fig:fastenerTrainingResults} shows the training error for the various ML models considered as a function of the reduced basis dimension, $K$, and the number of model parameters, $N_w$. We observe that SPSD-LLS results in the highest training errors and that increasing the basis dimension no longer results in lower errors after around $\reducedDimensionCombined = 8$. In training, the standard LLS model is much more accurate than SPSD-LLS and results in accurate models that achieve a relative error around $10^{-3}$. We observe that NN is more accurate than LLS for low reduced basis dimensions, but is surprisingly \textit{less accurate} for larger basis dimensions. Given that the data are inherently nonlinear due to contact, we theorize that this result is due to LLS over-fitting to the data when there are many model parameters; note that LLS is analytically solved and is not subject to any optimization process like NN. For small reduced basis dimensions, where the models have fewer parameters, NN presumably outperforms LLS due to its ability to model nonlinear responses. Lastly, we observe that SPSD-NN is significantly more accurate than SPSD-LLS and is competitive with both LLS and NN. This result demonstrates that we are able to generate a structure-preserving SPSD stiffness matrix model with a comparable accuracy to LLS and NN.  It is interesting to observe that the improvement between SPSD-NN and SPSD-LLS is far greater than between NN and LLS.

\subsection{Coupled FEM-ML results}
\begin{figure}
\begin{center}
\begin{subfigure}[t]{0.48\textwidth}
\includegraphics[trim={0cm 0cm 0cm 0cm},clip,width=1.0\linewidth]{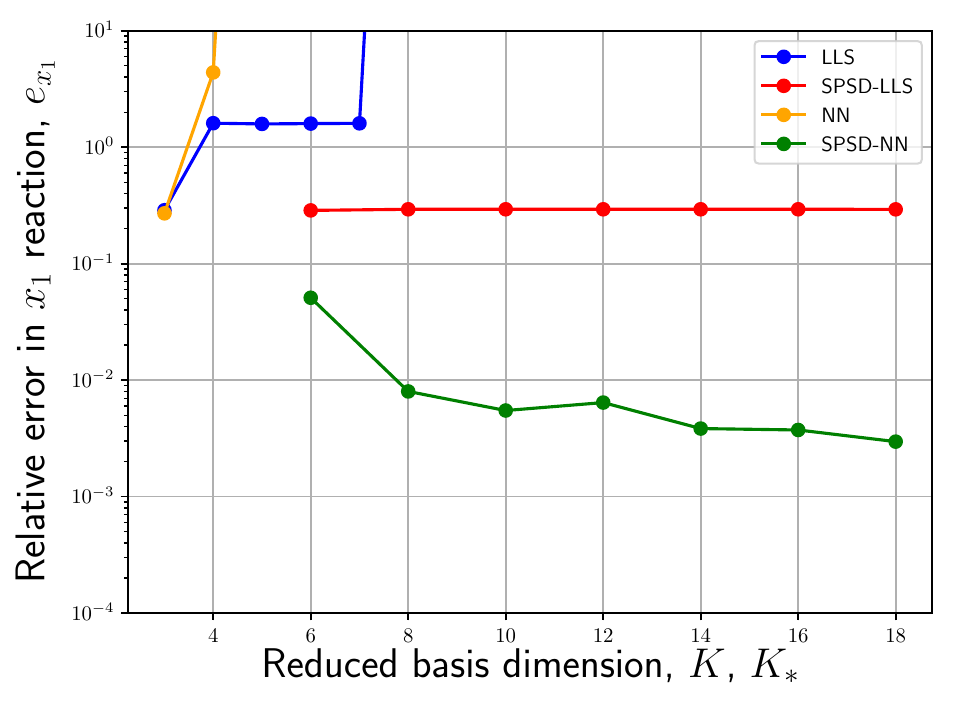}
\end{subfigure}
\begin{subfigure}[t]{0.48\textwidth}
\includegraphics[trim={0cm 0cm 0cm 0cm},clip,width=1.0\linewidth]{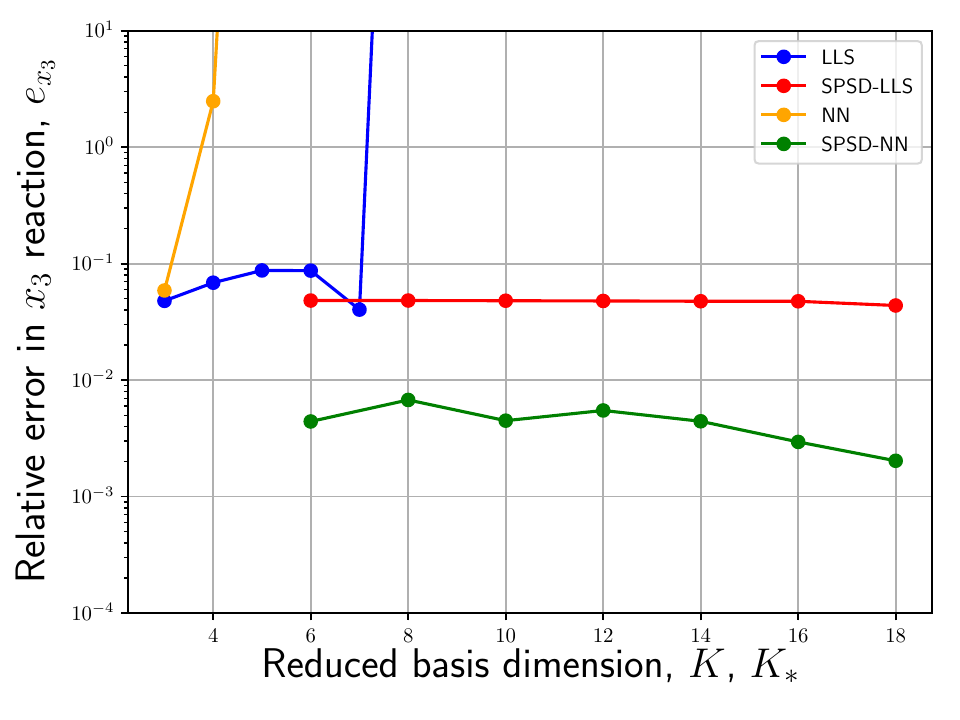}
\end{subfigure}
\caption{\contactCaption\ Coupled FEM-ML results for various ML models as a function of reduced basis dimension. Note that the SPSD models employ the combined basis of dimension $K^* = 2K$.}
\label{fig:fastenerResultOne}
\end{center}
\end{figure}

We now examine results for when the ML model is coupled to the FEM solver. We consider 9 runs for $$\beta = \{-0.005, -0.00475,-0.0045,-0.00425,-0.0025,-0.00125,-0.001,-0.00025,0.000\}.$$ This set of runs contains configurations that were in the training set as well as novel testing configurations. Figure~\ref{fig:fastenerResultOne} shows the convergence of the relative error (summed over all cases) for the predicted $x_1$-reaction and $x_3$-reaction forces for the various models as a function of the reduced basis dimension. We observe that SPSD-NN is by far the best performing method. We observe that it remains stable and becomes more accurate as the model complexity grows; the most complex model tested led to relative errors of around $0.2\%$. We next observe that LLS and NN go unstable as the model complexity increases. Further, even when stable, both methods associate with high relative errors around $50\%$. Lastly, we observe that SPSD-LLS remains stable as the reduced basis dimension grows, but similarly to Figure~\ref{fig:fastenerTrainingResults}, the model is unable to improve as we allow it to have more parameters. This lack of improvement can be attributed to the fact that the data are nonlinear and the SPSD-LLS model cannot be overfit to the data in the same way that LLS can due to the structure we impose on the model. The SPSD-NN model does yield improved results as the reduced basis dimension grows. 

\begin{figure}
\begin{center}
\begin{subfigure}[t]{0.48\textwidth}
\includegraphics[trim={0cm 0cm 0cm 0cm},clip,width=1.0\linewidth]{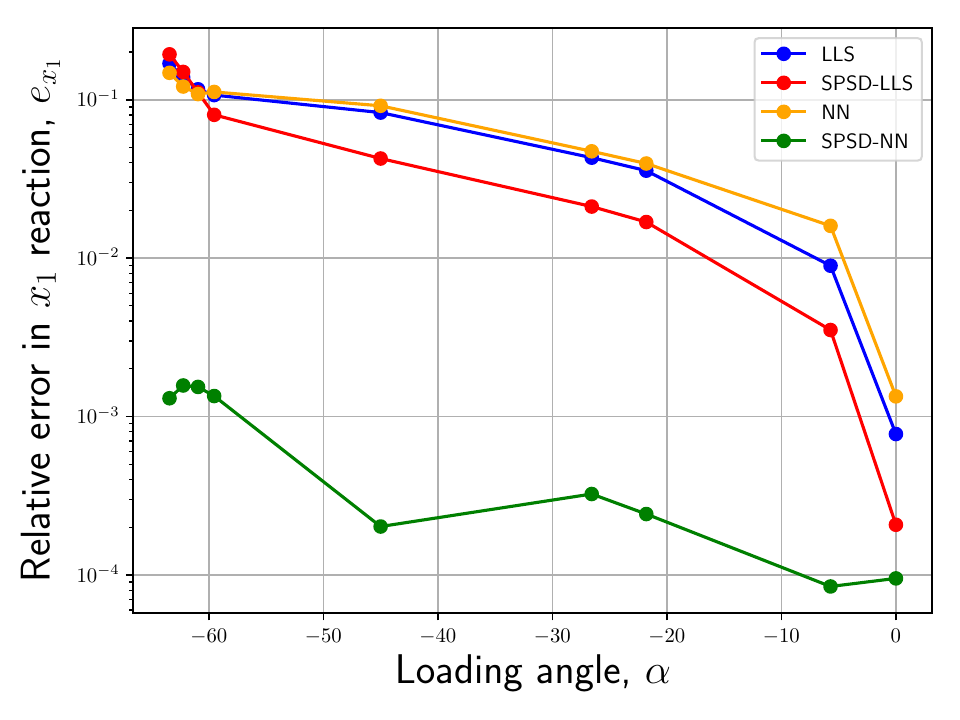}
\end{subfigure}
\begin{subfigure}[t]{0.48\textwidth}
\includegraphics[trim={0cm 0cm 0cm 0cm},clip,width=1.0\linewidth]{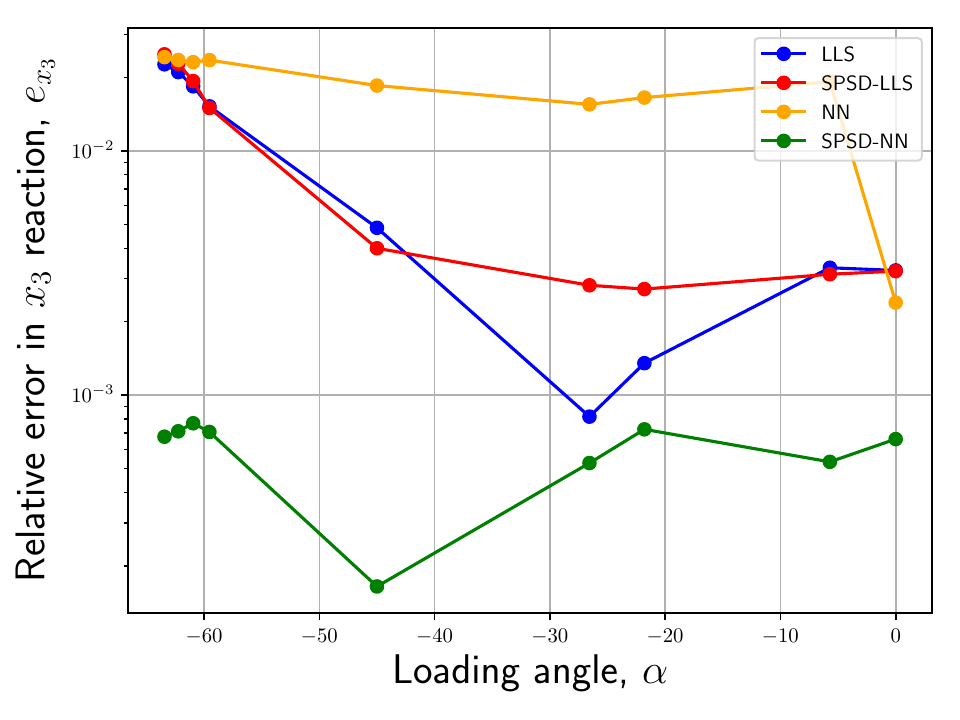}
\end{subfigure}
\caption{\contactCaption\ Coupled FEM-ML results for various ML models as a function of loading angle. Results are shown for the optimal configuration for each model, which corresponds to $K=3$, $K^*=6$, $K=3$, and $K^*=18$ for the LLS, SPSD-LLS, NN, and SPSD-NN models, respectively.}
\label{fig:fastenerResultTwo}
\end{center}
\end{figure}

Figure~\ref{fig:fastenerResultTwo} shows the relative error of the best performing model of each method as a function of loading angle for the various cases; the best performing models correspond to $K=3$, $K^*=6$, $K=3$, and $K^*=18$ for the LLS, SPSD-LLS, NN, and SPSD-NN models, respectively. We note that contact occurs in the FEM for all cases when, approximately, $|\alpha| > 40$. As expected \SPSDNN\ is by far the best performing method and leads to relative errors of under $0.2\%$ in all cases. \SPSDLLS\ is the next best performing method, but in general \SPSDLLS, LLS, and NN all associate with relatively high errors. This is particularly true when there is contact in which case all methods associate with $>10\%$ error when it comes to predicting the $x_1$-reaction force.  

\begin{figure}
\begin{center}
\begin{subfigure}[t]{1.0\textwidth}
\includegraphics[trim={0cm 0cm 0cm 0cm},clip,width=1.0\linewidth]{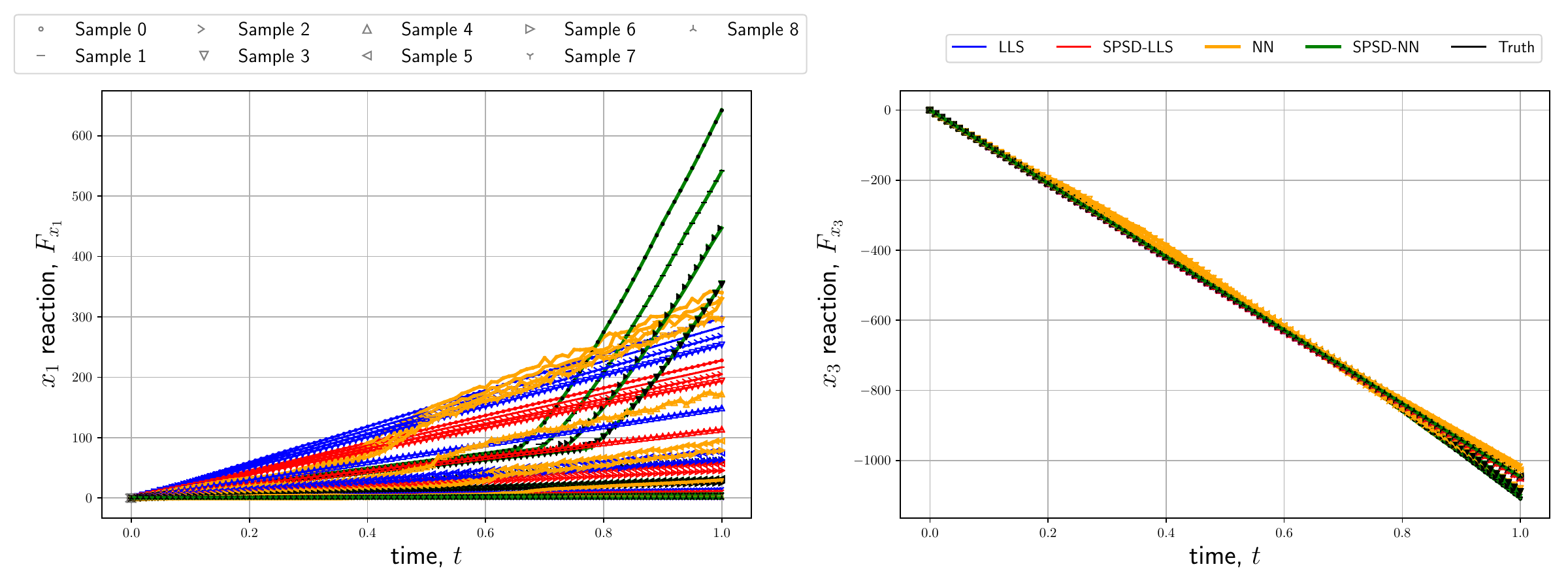}
\end{subfigure}
\caption{\contactCaption\ Coupled FEM-ML results for all testing configurations. We show the measured reaction force in the $x_1$ (shear) direction on the left and the measured reaction force in the $x_3$ (tensile) direction on the right. We note that contact is more predominant in the shear direction. Values with small absolute loading angles correspond to the cases with low reaction forces in the shear direction, while values with large absolute loading angles correspond to cases with large loading angles in the shear direction.}
\label{fig:fastenerResultAll}
\end{center}
\end{figure}

Figure~\ref{fig:fastenerResultAll} shows the reaction forces predicted by the coupled FEM-ML models for all configurations. 
We observe that \SPSDNN\ is, by far, the best performing method and leads to predictions that lie on top of the truth values. \SPSDNN\ is able to correctly capture contact and correctly predict the change in the stiffness. For larger reaction forces, which correspond to larger (absolute) loading angles, \SPSDNN\ is able to capture the ramp-like behavior of the reaction force. The other nonlinear method, NN, is unable to characterize this behavior despite it performing similar to \SPSDNN\ in training. For lower absolute loading angles, in which case contact does not occur and the response is mostly linear, \SPSDNN\ is able to still predict the correct stiffness without being ``contaminated" by its ability to capture the increased stiffness once contact occurs. All other models are unable to capture this behavior. Both LLS and \SPSDLLS, which are linear models, are unable to model the contact non-linearity and simply bisect its behavior. This leads to an under-prediction of the stiffness in the post-contact regime and an over-prediction of the stiffness in the pre-contact regime. The standard NN model attempts to capture the contact but clearly is subject to instabilities throughout the solve. 

\begin{figure}
\begin{center}
\begin{subfigure}[t]{0.75\textwidth}
\includegraphics[trim={3cm 1cm 5cm 1cm},clip,width=1.0\linewidth]{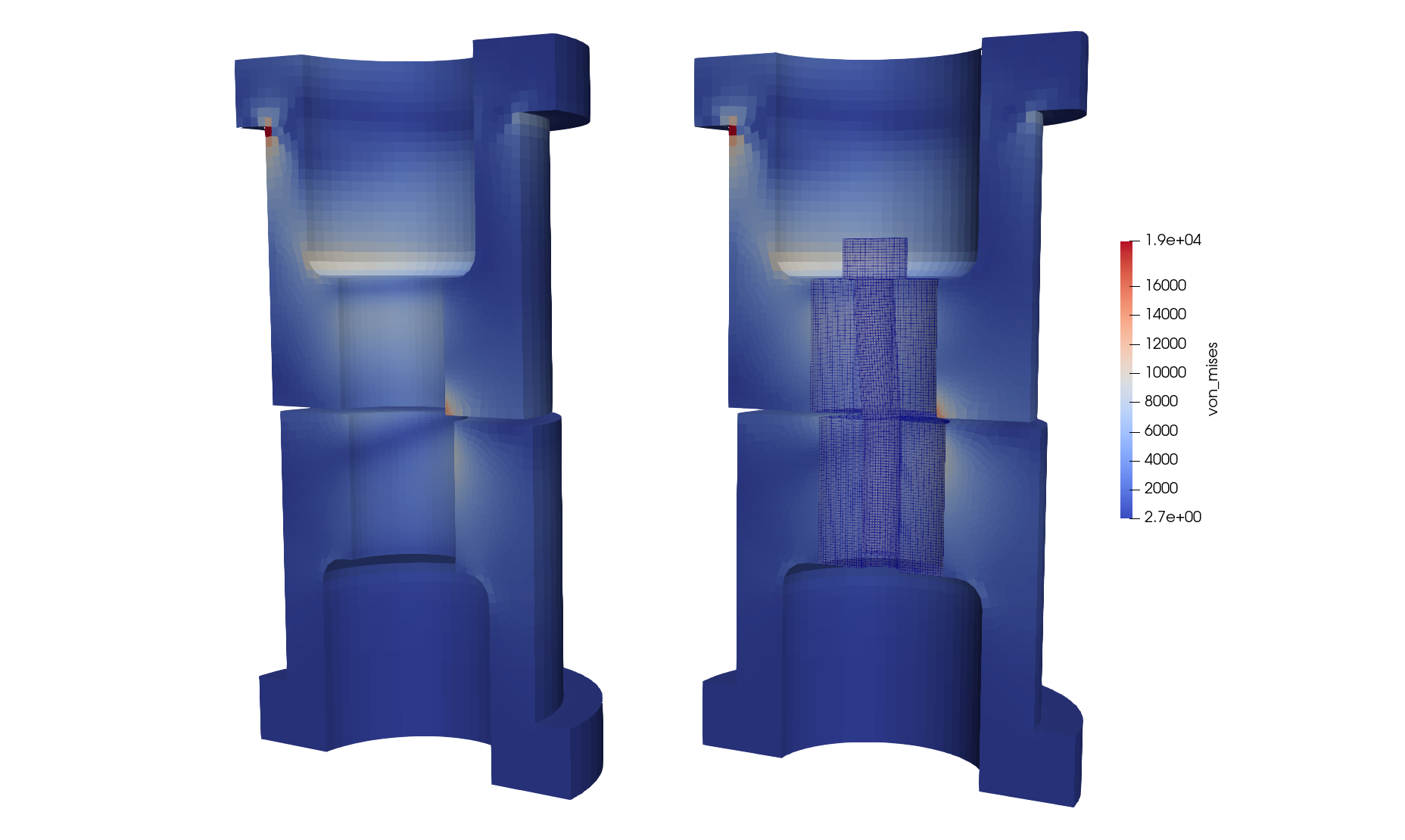}
\end{subfigure}
\caption{\contactCaption\ Predictions for the max von Mises stresses by the FEM-ML (left) and FEM-only (right) models. Note that the displacements are enlarged by a factor of 10 for visualization.}
\label{fig:fastenerResultFour}
\end{center}
\end{figure}

\begin{figure}
\begin{center}
\begin{subfigure}[t]{0.455\textwidth}
\includegraphics[trim={0 0 0 0},clip,width=1.0\linewidth]{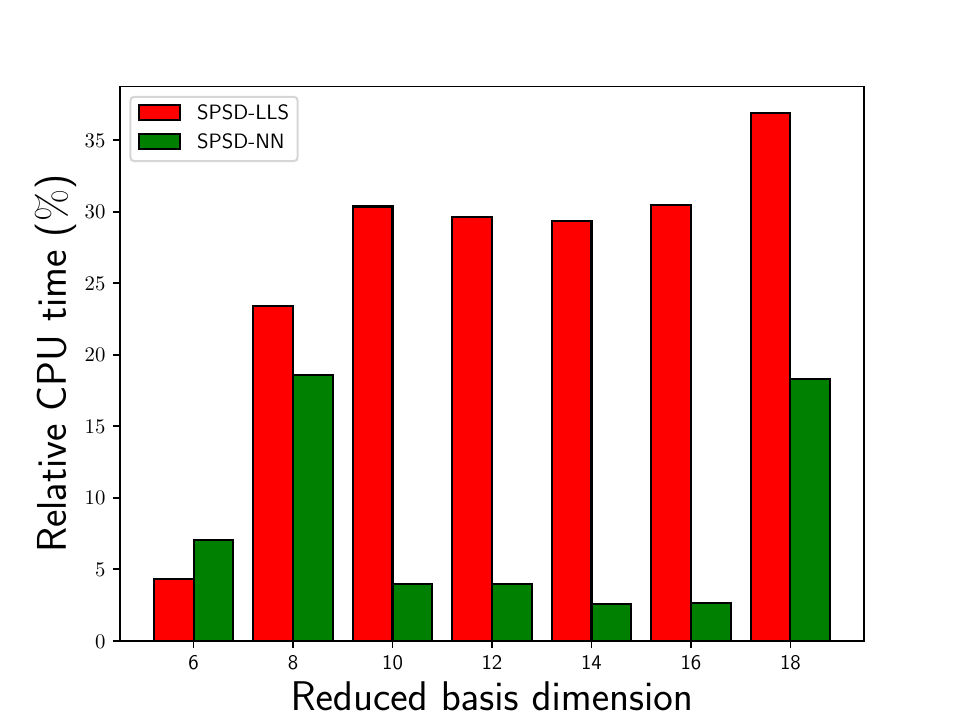}
\end{subfigure}
\caption{\contactCaption\ Relative CPU times for the FEM--ML coupled reduced-order models (i.e., CPU time of the FEM-ML model divided by CPU time of the FEM-only full-order model which resolves the fastener). The reported CPU times are averaged across the nine testing configurations.} 
\label{fig:fastenerResultTimes}
\end{center}
\end{figure}

Figure~\ref{fig:fastenerResultFour} shows contour plots of the max von Mises stress throughout the outer bushing as predicted by the coupled FEM-ML model with the best performing \SPSDNN\ model (left) and FEM-only model (right). For this range of von Mises stress contours, the coupled FEM-ML model leads to results that are visually identical to the FEM-only model. Lastly, Figure~\ref{fig:fastenerResultTimes} reports the relative CPU times of the FEM--ML reduced-order models as compared to the FEM--only full-order model. We note that we only show CPU times for the \SPSDLLS\ and \SPSDNN\ models given that LLS and NN are not stable for most configurations. We observe that the \SPSDLLS\ model results in speedups between 3-20x, depending on the model configuration, while \SPSDNN\ results in speedups between 5-25x, again depending on the model configuration. Interestingly, we observe that the \SPSDNN\ model is, on average, faster than the \SPSDLLS\ model. An investigation into this reveals that, in the conjugate gradient solver, converging the FEM--ML model with the \SPSDNN\ model requires, on average, significantly fewer iterations than the \SPSDLLS\ model. This result suggests that the \SPSDNN\ model is producing better conditioned stiffness matrices, and highlights adding penalties on the conditioning of the stiffness matrix as an interesting path for future work. Lastly, we highlight that we observe a slight increase in CPU time for \SPSDLLS\ as the reduced basis dimension is increased, suggesting that conditioning of the resulting stiffness matrix degrades as the basis dimension is increased. No such pattern is apparent for the \SPSDNN\ model. 

\subsection{Fastener undergoing three-dimensional radial loading with contact and preload}
\begin{figure}
\begin{center}
\begin{subfigure}[t]{0.8\textwidth}
\includegraphics[trim={6cm 4cm 4cm 6cm},clip,width=1.0\linewidth]{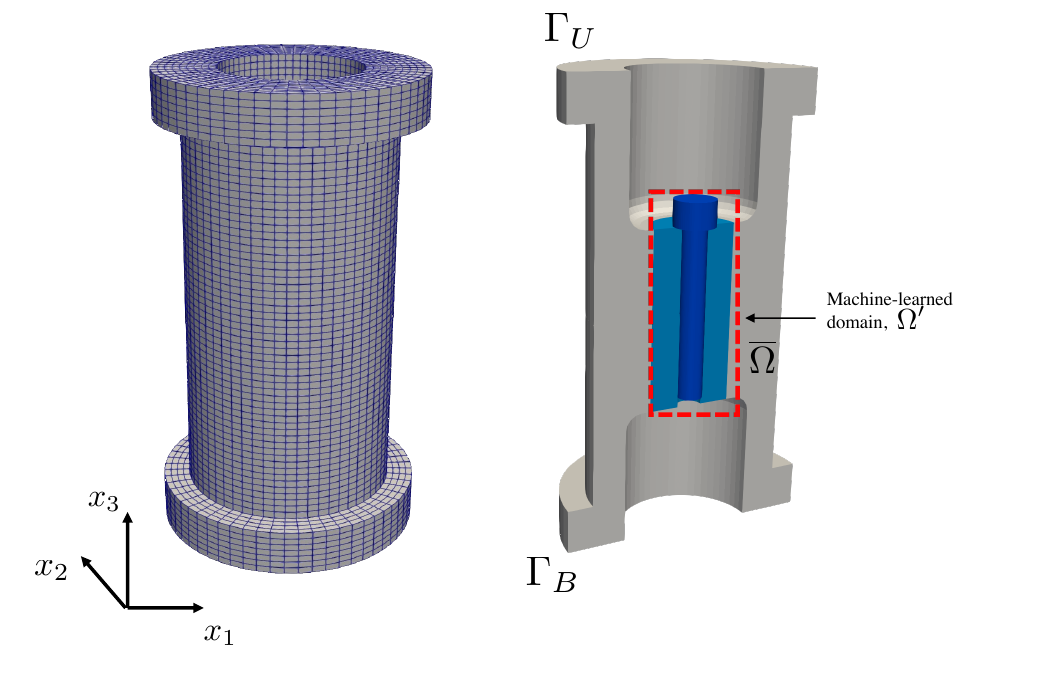}
\end{subfigure}
\caption{\preloadCaption\ Schematic of the problem configuration. On the left we show the outer computational mesh and on the right we show the geometry clipped down the middle. The middle region highlighted by the red box comprises the fastener-bushing domain that will be replaced with the ML model. }
\label{fig:preloadSchematic}
\end{center}
\end{figure}

The final problem we consider is again a fastener-bushing geometry, but this time we consider a three-dimensional loading profile with an initial preloaded state. The problem configuration is shown in Figure~\ref{fig:preloadSchematic}. The full FEM mesh comprises $N = 234,750$ degrees of freedom. We remove the immediate domain around the fastener, which results in an interface with $\nCoarseNodesDomainBoundaryOneTwo = 6840$ nodes. We again emphasize that contact occurs exclusively within the removed domain, and as a result our ML model must be able to accurately characterize contact. 
We again employ homogeneous Dirichlet conditions on the bottom boundary. The quasi-static loading profile for $\x$ on $\Gamma_U$ is given, in inches, by 
\begin{equation}
\begin{split}
&{\displacementFEM}_{1}(t,\x) = 
\beta t  \qquad \\
&{\displacementFEM}_{2}(t,\x) = 
\alpha t   \\
&{\displacementFEM}_{3}(t,\x) = 
\frac{1}{400}  t. \\
\end{split}
\end{equation}
We consider 121 training configurations for $\beta = \{-0.005 + 0.001 i\}_{i=0}^{10}$, $\alpha = \{-0.005 + 0.001 i\}_{i=0}^{10}$, and $t \in \{0.01i\}_{i=1}^{101}$ $\mathsf{s}$. This dataset comprises quasi-static trajectories with 101 quasi-static time steps per trajectory for pulls ranging from approximately $-65^{\circ}$ to $65^{\circ}$ in both directions. In addition, for this example we consider a preloaded state such that, at time $t=0$, the initial force on the interface is $\mathbf{f}_{\mathsf{interface}} = \mathbf{f}_{\mathsf{preload}} \ne 0$. This preload is applied through the use of an artificial strain in the axial direction of the fastener. 

We employ the same material configuration as the previous two exemplars. For $\x \in \domainOne$ the material model is characterized by a Young's modulus of 28.5e6 \textsf{psi} with Poisson's ratio of $\nu = 0.3$. For $\x \in \domainTwo$ the material model is characterized by a Young's modulus of $29 \times 10^6$ \textsf{psi} with a Poisson's ratio of $\nu=0.3$. As a QoI, we again consider the integrated reaction forces on the exterior of the bottom bushing, but this time in the $x_1$, $x_2$, and $x_3$ directions.  

For this complex exemplar we only investigate the performance of the \SPSDNN\ model given its superior performance over other formulations. We additionally note that we employ a slightly different training configuration for this example: the batch size is 60 and our initial learning rate is $\mathsf{lr} = 2.5 \times 10^{-4}$. The rest is the same as is described in Section~\ref{sec:ml_training}. 


\subsection{Offline training results}
\begin{figure}
\begin{center}
\begin{subfigure}[t]{0.49\textwidth}
\includegraphics[trim={0cm 0cm 0cm 0cm},clip,width=1.0\linewidth]{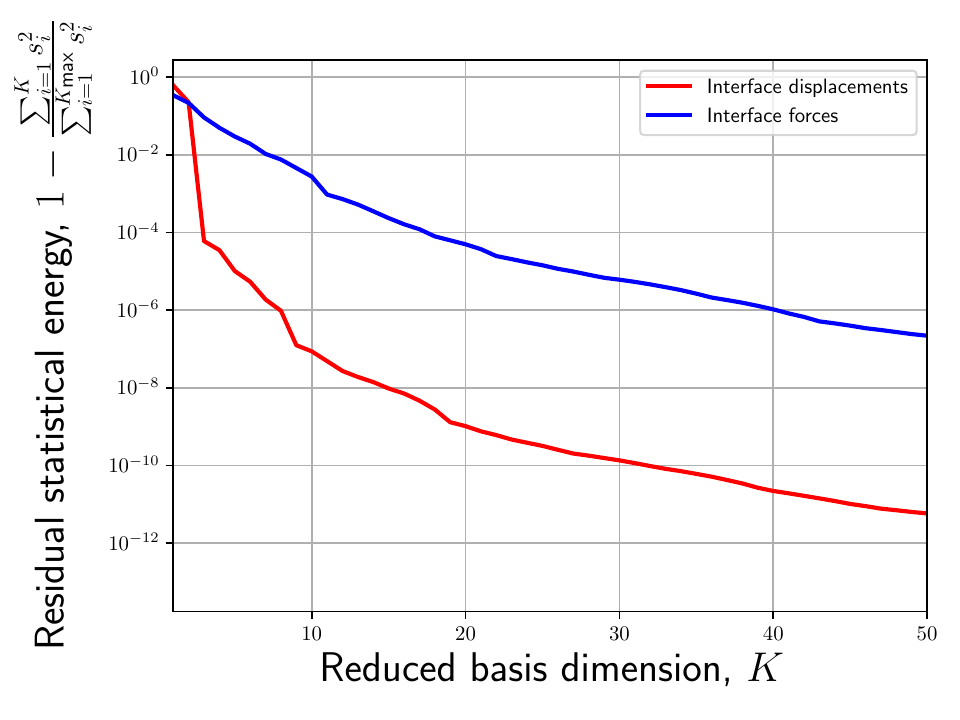}
\end{subfigure}
\caption{\preloadCaption\ Residual statistical energy of the reduced basis approximation for the interface displacement and force fields.}
\label{fig:preloadTrainingBasisConvergence}
\end{center}
\end{figure}

We again first assess the reduced basis approximations. To capture preload, the force offset vector is set to be the force state resulting from preload, i.e., $\forceOffset = \mathbf{f}_{\mathsf{preload}}$\footnote{The preload force vector is non-uniform and is extracted from the FEM simulation in the training phase. In future work, we plan to parameterize the preload force to handle different types of preload.}.  Figure~\ref{fig:preloadTrainingBasisConvergence} depicts the residual statistical energy associated with the reduced basis approximation. As compared to the previous examples, we observe a much slower decay in the residual statistical energy of the interface force field: 40 basis vectors are required to reach a residual statistical energy of $10^{-6}$, while in the previous two exemplars this tolerance was achieved with only 10 basis vectors. This result demonstrates the more sophisticated physics present in this example. Despite the slower decay in singular values associated with the interface forces, we are still able to identify a subspace of a much reduced dimension capable of representing the majority of the statistical energy of the force and displacement interface fields.

\subsection{Coupled FEM-ML results}
\begin{figure}
\begin{center}
\begin{subfigure}[t]{1.0\textwidth}
\includegraphics[trim={0 0 0 0},clip,width=1.0\linewidth]{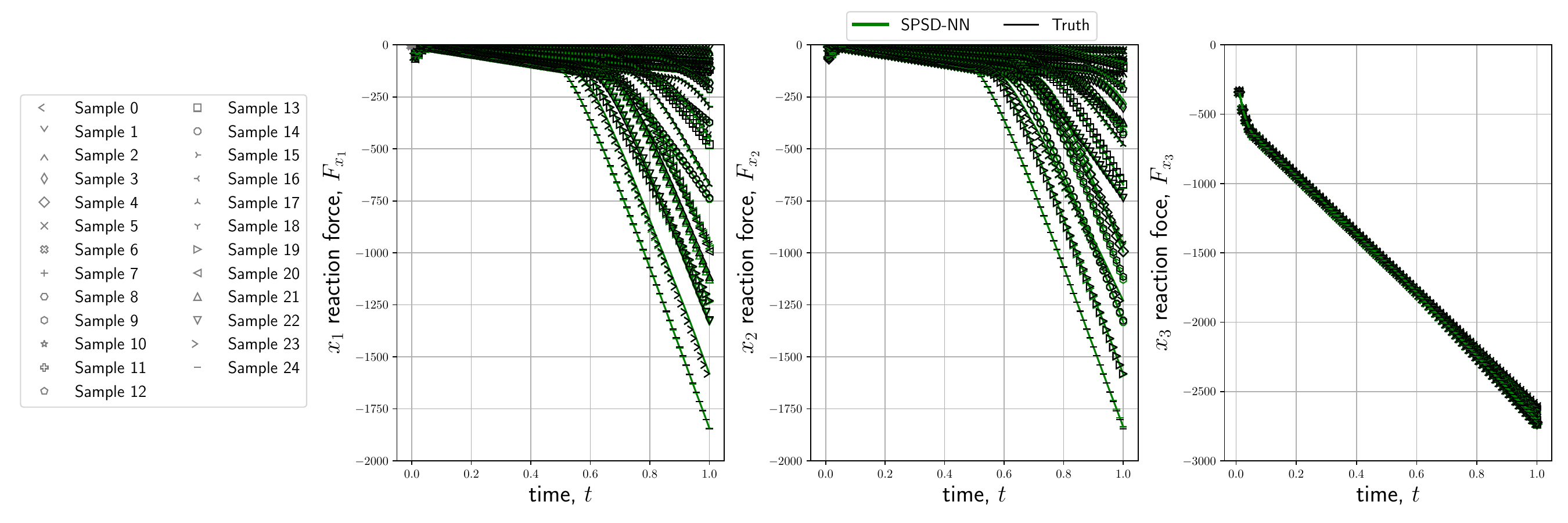}
\end{subfigure}
\caption{\preloadCaption\ Coupled FEM-ML results for all testing configurations. From left to right, we show the predicted reaction force in the $x_1$ (shear), $x_2$ (shear), and $x_3$ (tension) directions, respectively. Each line corresponds to a different quasi-static trajectory obtained from a different value of $\beta$ and $\alpha$.}
\label{fig:preloadOnlineResultOne}
\end{center}
\end{figure}

\begin{figure}
\begin{center}
\begin{subfigure}[t]{0.32\textwidth}
\includegraphics[trim={0 0 0 0},clip,width=1.0\linewidth]{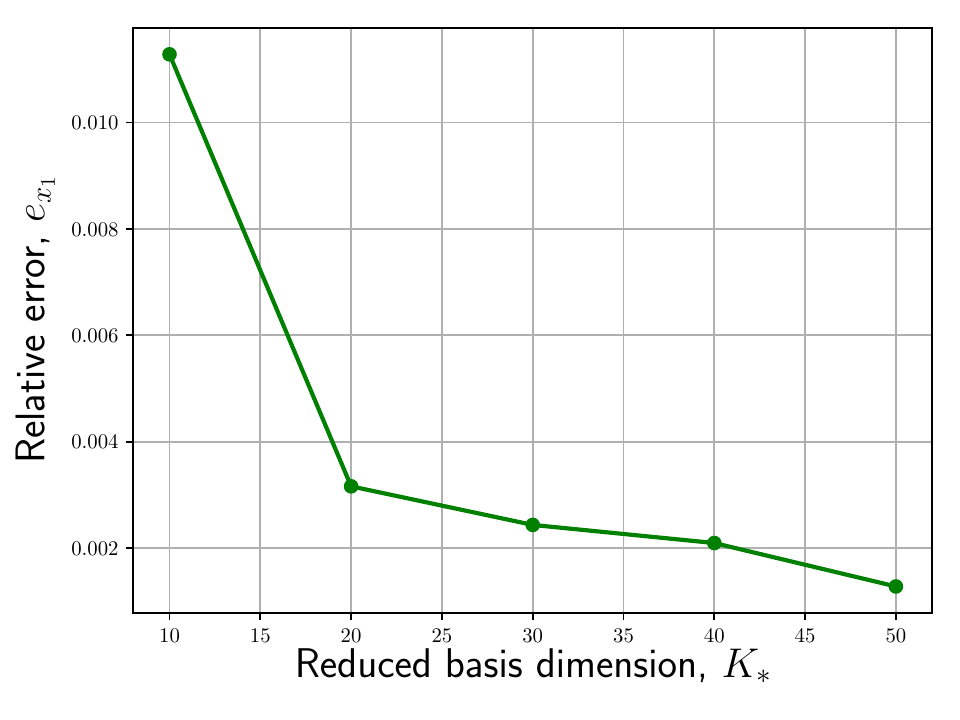}
\end{subfigure}
\begin{subfigure}[t]{0.32\textwidth}
\includegraphics[trim={0 0 0 0},clip,width=1.0\linewidth]{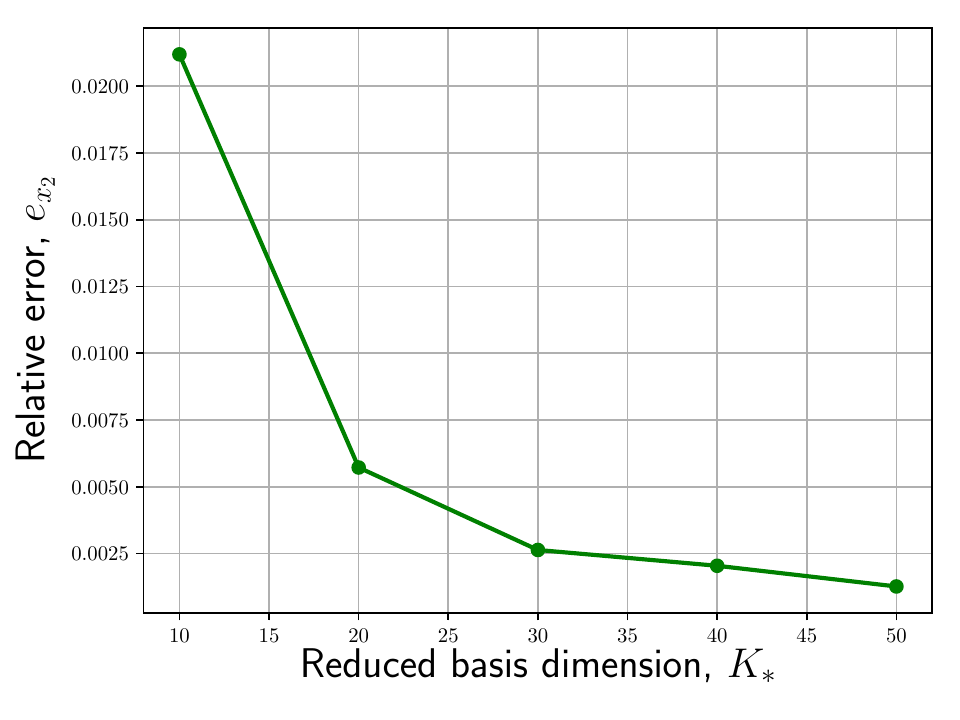}
\end{subfigure}
\begin{subfigure}[t]{0.32\textwidth}
\includegraphics[trim={0 0 0 0},clip,width=1.0\linewidth]{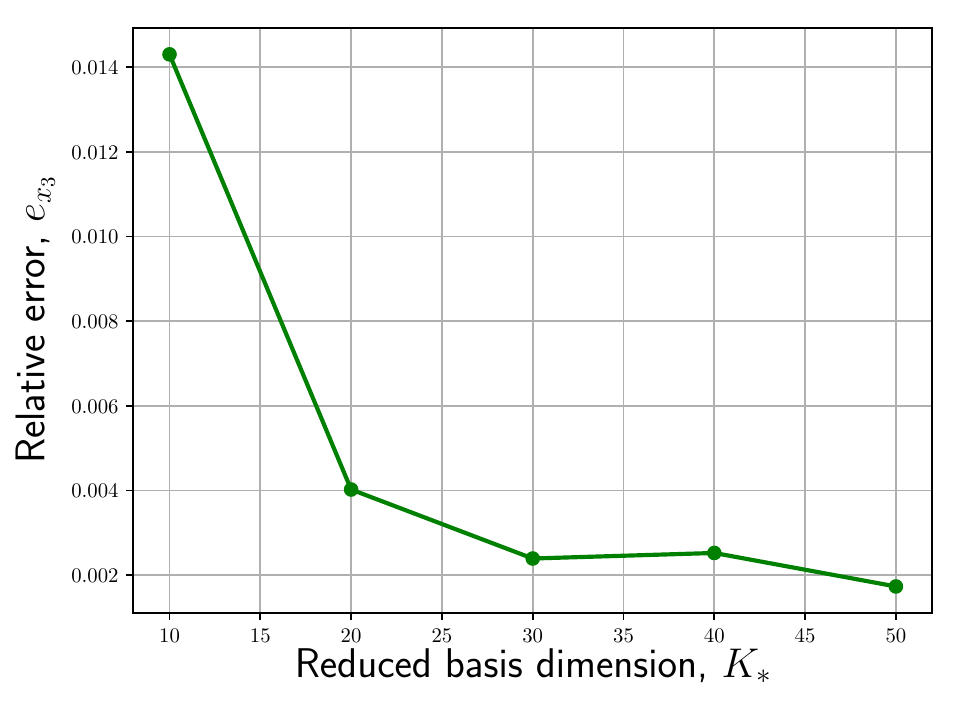}
\end{subfigure}
\caption{\preloadCaption\ Relative errors of coupled FEM-ML results with the \SPSDNN\ ML model as a function of reduced basis dimension. From left to right, we show the predicted reaction force in the $x_1$ (shear), $x_2$ (shear), and $x_3$ (tension) directions, respectively.}
\label{fig:preloadOnlineResultConvergence}
\end{center}
\end{figure}

We move directly to the coupled FEM--ML results. We consider a set of out-of-sample testing runs for $\beta = \{0.0005 + 0.001 i\}_{i=0}^{4}$, $\alpha = \{0.0005 + 0.001 i\}_{i=0}^{4}$. Figure~\ref{fig:preloadOnlineResultOne} shows QoI predictions of $K^*=50$ across all testing samples. Focusing on the reaction force in the $x_3$ direction, we observe that the ML model is able to replicate the initial preloaded state by virtue of the constant offset vector. We further observe that the ML model is able to simulate the change in stiffness that occurs due to the initial preload around $t \approx 0.05$. Next, examining the reaction forces in the $x_1$ and $x_2$ directions, we observe the ML model is able to accurately characterize the change in stiffness due to contact. The model correctly predicts both when contact occurs, as well as the magnitude of the resulting reaction force. We emphasize that all results are for out-of-sample configurations.  Next, Figure~\ref{fig:preloadOnlineResultConvergence} shows the convergence of the reaction QoI errors as a function of the reduced basis dimension. We observe monotonic convergence in both the $x_1$ and $x_2$ reaction errors, and almost monotonic convergence in the $x_3$ reaction error. 

\begin{figure}
\begin{center}
\begin{subfigure}[t]{0.455\textwidth}
\includegraphics[trim={0 0 0 0},clip,width=1.0\linewidth]{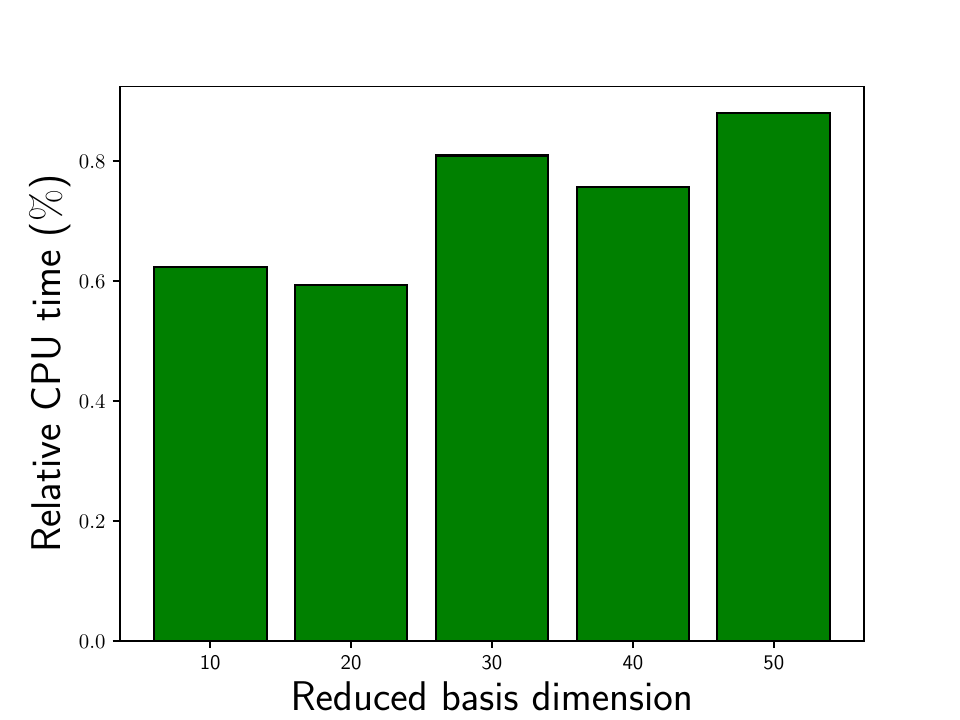}
\end{subfigure}
\caption{\preloadCaption\ 
Relative CPU times for the FEM--ML coupled reduced-order models (i.e., CPU time of the FEM-ML model divided by CPU time of the FEM-only full-order model which resolves the fastener). The reported CPU times are averaged across all testing configurations.}
\label{fig:preloadOnlineResultTimes}
\end{center}
\end{figure}

\begin{figure}
\begin{center}
\begin{subfigure}[t]{0.85\textwidth}
\includegraphics[trim={5cm 2cm 10cm 2cm},clip,width=1.0\linewidth]{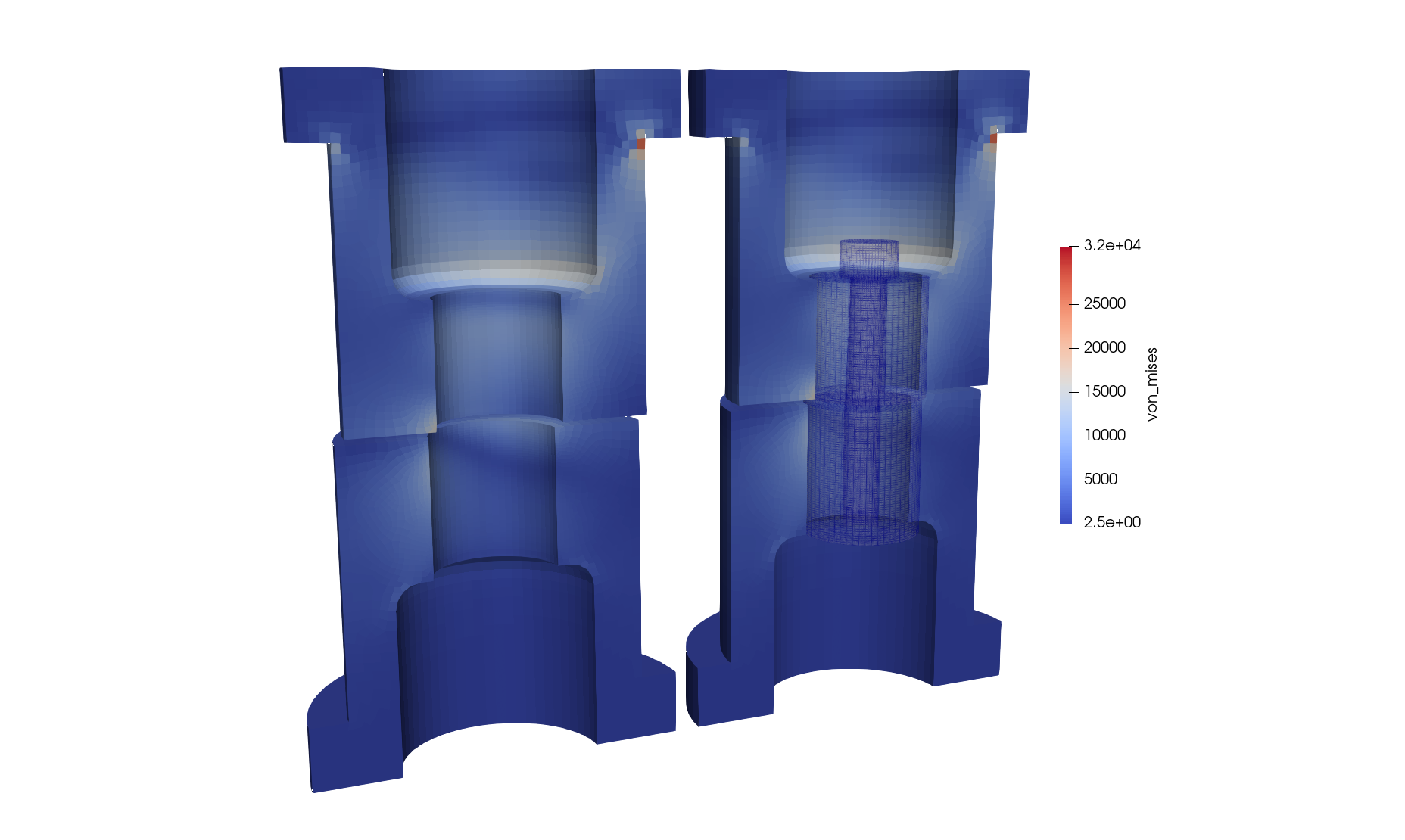}
\end{subfigure}
\caption{\preloadCaption\ Predictions for the max Von Mises stresses by the FEM-ML (left) and FEM-only (right) models for $\beta = 0.0045$, $\alpha = 0.0035$. Note that the displacements are enlarged by a factor of 10 for visualization. The FEM--ML solution was obtained approximately 130x faster than the FEM-only solution.}
\label{fig:fastenerPreloadContour}
\end{center}
\end{figure}

Figure~\ref{fig:preloadOnlineResultTimes} presents the relative wall-times of the various FEM--ML coupled models as compared to the FEM-only model. We observe that, for all basis dimensions, the FEM--ML coupled systems result in over 100x speedups. We again do not observe a noticeable change in cost with increasing the reduced basis dimension. Figure~\ref{fig:fastenerPreloadContour} shows contour plots for the max von Mises stress as predicted by the FEM--ML, with the \SPSDNN\ model at a reduced basis dimension of $\reducedDimensionCombined=40$, and FEM-only models. No differences in results are distinguishable.  

\begin{figure}
\begin{center}
\begin{subfigure}[t]{0.32\textwidth}
\includegraphics[trim={0 0 0 0},clip,width=1.0\linewidth]{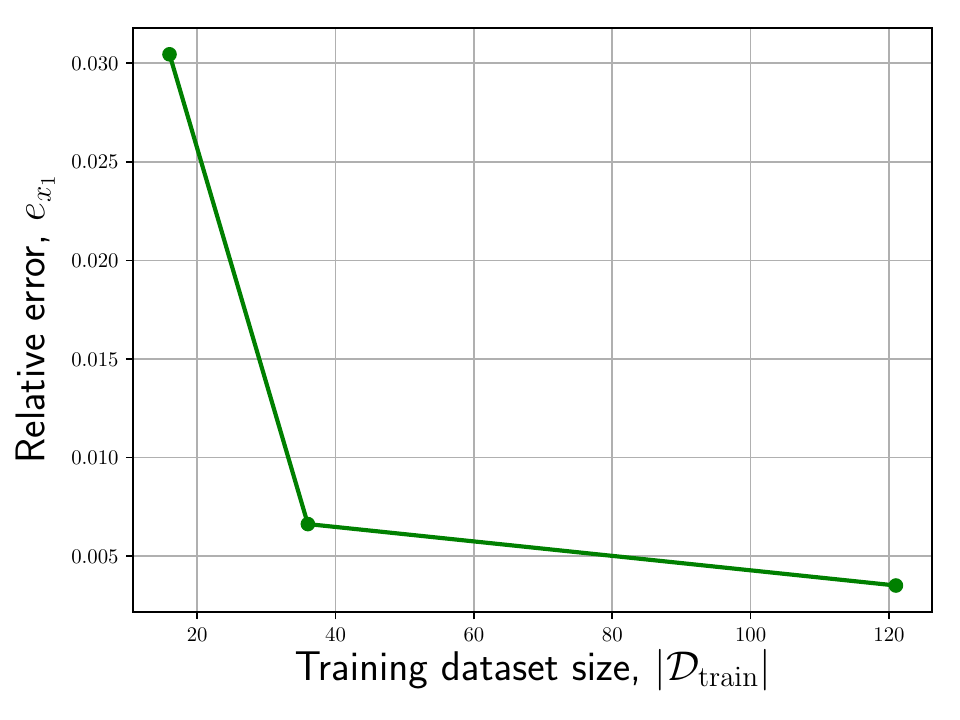}
\end{subfigure}
\begin{subfigure}[t]{0.32\textwidth}
\includegraphics[trim={0 0 0 0},clip,width=1.0\linewidth]{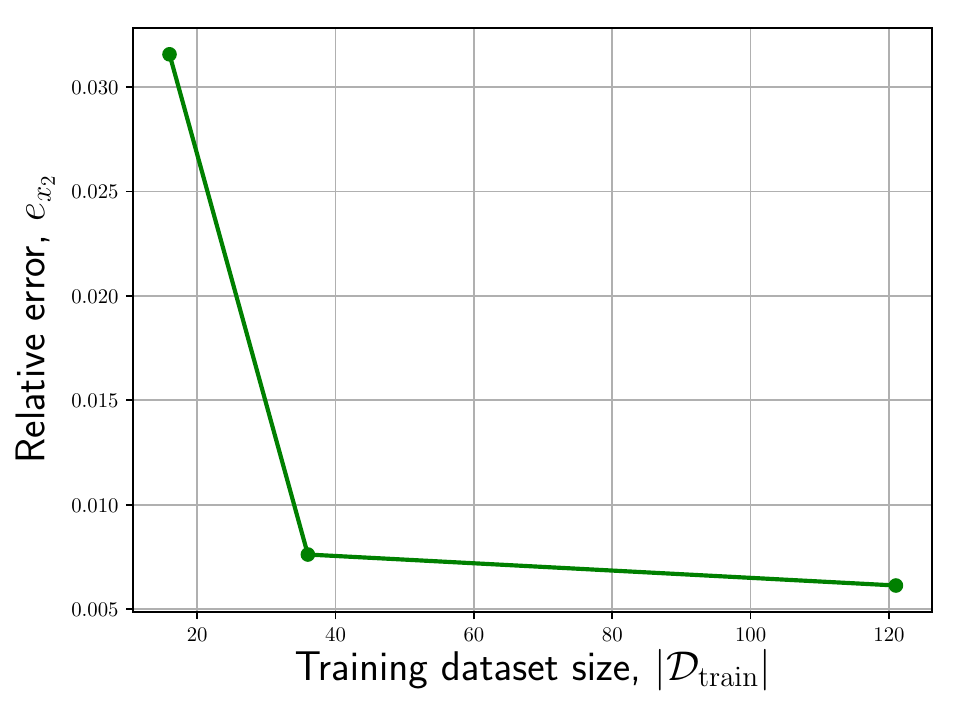}
\end{subfigure}
\begin{subfigure}[t]{0.32\textwidth}
\includegraphics[trim={0 0 0 0},clip,width=1.0\linewidth]{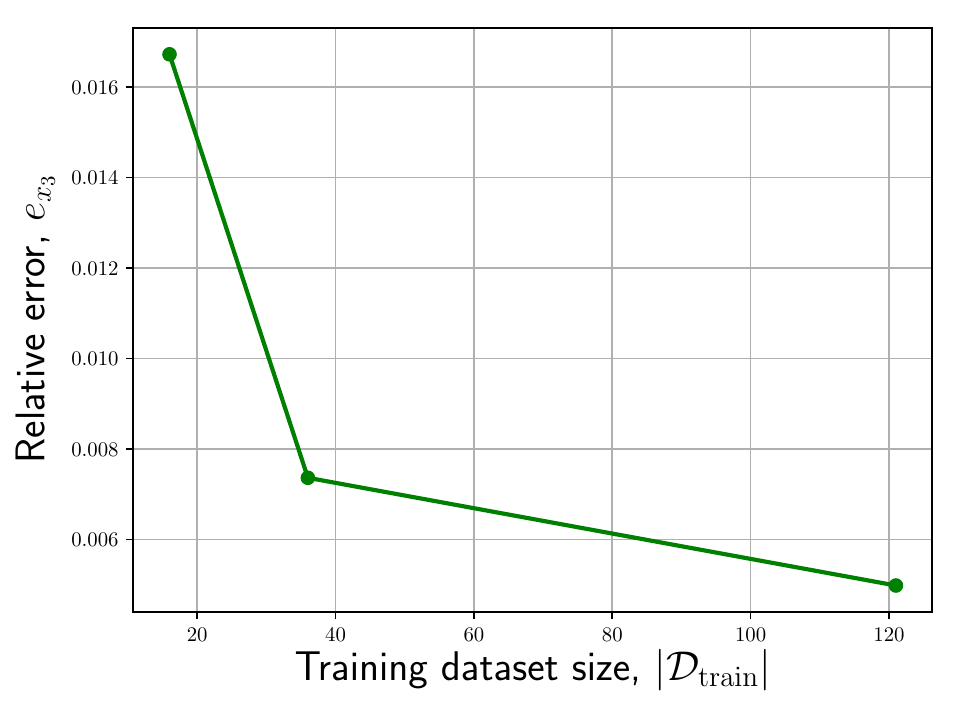}
\end{subfigure}
\caption{\preloadCaption\ Relative errors of coupled FEM-ML results with the \SPSDNN\ ML model as a function of training data set size. From left to right, we show the predicted reaction force in the $x_1$ (shear), $x_2$ (shear), and $x_3$ (tension) directions, respectively.}
\label{fig:preloadOnlineResultConvergenceData}
\end{center}
\end{figure}

Lastly we consider the impact of the training dataset size on model performance. To this end, we consider ML models with a reduced basis dimension of $\reducedDimensionCombined=20$ trained using three different training sets, 
\begin{enumerate}
\item  $\beta = \{-0.005 + 0.0005 i\}_{i=0}^{10}$, $\alpha = \{-0.005 + 0.005 i\}_{i=0}^{10}$,
\item  $\beta = \{-0.005 + 0.0005 i\}_{i=0,2,4,6,8,10}$, $\alpha = \{-0.005 + 0.005 i\}_{i=0,2,4,6,8,10}$,
\item  $\beta = \{-0.005 + 0.0005 i\}_{i=0,3,6,9}$, $\alpha = \{-0.005 + 0.005 i\}_{i=0,3,6,9}$.
\end{enumerate}
All three data sets employ $t \in \{0.01i\}_{i=1}^{101}$ $\mathsf{s}$. For each case we test on the same out-of-sample testing runs as before. Figure~\ref{fig:preloadOnlineResultConvergenceData} shows the convergence of the reaction QoI errors as a function of the training dataset size, denoted by $|\mathcal{D}_{\mathrm{train}}|$. We observe monotonic convergence in all QoIs as the number of training data increase. For the coarsest training dataset, we observe relatively large QoI errors around $3\%$, while for the finest dataset we observe sub $1\%$ QoI errors. As expected, these results indicate that enriching the training set size improves model accuracy. 

\section{Conclusion}\label{sec:conclude}
This work introduced a machine-learning strategy for finite element analysis of solid mechanics wherein we replaced hard-to-resolve portions of a computational domain with a data-driven surrogate. We proposed two types of data-driven surrogates: one that maps directly from the interface displacements to the interface forces, and one that maps from the interface displacements to the interface forces by virtue of a stiffness matrix that is enforced to be SPSD. We demonstrated, in a simplified setting, that this latter formulation results in a global coarse-scale problem that is symmetric positive-definite (SPD) which guarantees a unique solution and makes the form amenable to conjugate gradient solvers.

We presented numerical experiments across three exemplars spanning a range of physics. These exemplars demonstrated that, in general, direct force-to-displacement models are not robust when combined with a conjugate gradient solver. At this time it is not clear if this lack of robustness is due to interplay between the model form and the conjugate gradient solver, or if rather it is due to the resulting problem being ill-posed. Even if model robustness can be attributed to the conjugate gradient solver, in which case improved performance could be obtained, e.g., by switching to a general minimized residual-based solver, this is still a significant issue given the wide use and effectiveness of conjugate gradient-based solvers in SM.

Our numerical experiments demonstrated that the displacement-to-force via an SPSD stiffness matrix approach resulted in robust and efficient models. In all of our testing cases, the FEM-ML model with the \SPSDNN\ approach was able to predict QoIs with a relative error of less than $0.5\%$. On the final exemplar, which was the most complex test case, this formulation led to sub $0.25\%$ errors with speed-ups of over 100x as compared to the standard FEM simulation.

The results of this work demonstrate the promise of using machine-learning for alleviating the computational burden associated with hard-to-resolve portions of a computational domain, like threaded fasteners, for traditional finite element methods. Follow on work will focus on aspects including (1) extension to more sophisticated material models, e.g., elastic-plastic, where history behavior may become important, (2) extension to component-level exemplars with fasteners, where the same fastener model is repeated throughout the domain, (3) adding additional physical constraints to our models, e.g., rigid body rotation, and (4) prediction of QoIs within the substituted subdomain, such as fastener stress and fastener failure.

\section{Acknowledgements}
This paper describes objective technical results and analysis. 
The authors acknowledge support from Sandia Advanced Simulation and Computing projects 65755 and 103723.  
Any subjective views or opinions that might be expressed in the paper do
not necessarily represent the views of the U.S. Department of Energy
or the United States Government. Sandia National Laboratories is a
multimission laboratory managed and operated by National Technology \&
Engineering Solutions of Sandia, LLC, a wholly owned subsidiary of
Honeywell International Inc., for the U.S. Department of Energy’s
National Nuclear Security Administration under contract DE-NA0003525.

\bibliographystyle{siam}
\bibliography{refs}

\newpage
\begin{appendix}

\section{PyTorch Code}\label{sec:pytorch_code}
This appendix provides the PyTorch code used for defining the neural network and symmetric positive definite neural network models. Listing~\ref{listing_one} provides the PyTorch code for the vanilla neural network architecture, while Listing~\ref{listing_two} provides code for the SPSD neural network architecture. 

\begin{lstlisting}[language=Python,caption=Code for neural network architecture used for the Neural Network surrogate model.,label=listing_one]
import numpy as np
import torch
import torch.nn.functional as F
import torch.nn as nn
torch.set_default_dtype(torch.float64)

class NN(nn.Module):
    def __init__(self,nHiddenLayers,nNeuronsPerLayer,nFeatures,nOutputs):
        super(NN, self).__init__()
        forward_list = []
        self.numHiddenLayers = nHiddenLayers
        self.numLayers = self.numHiddenLayers + 1
        dim = np.zeros(nHiddenLayers+2,dtype='int')
        dim[0] = nFeatures
        for i in range(1,nHiddenLayers+1):
         dim[i] = nNeuronsPerLayer
        dim[-1] = nOutputs
        self.dim = dim
        input_dim = dim[0:-1]
        output_dim = dim[1::]

        for i in range(0,nHiddenLayers+1):
          forward_list.append(nn.Linear(input_dim[i], output_dim[i]))

        self.forward_list = nn.ModuleList(forward_list)
        self.activation = F.relu

    def forward(self,x):
      for i in range(0,self.numLayers-1):
        x = self.activation(self.forward_list[i](x))
      x = self.forward_list[-1](x)
      return x
\end{lstlisting}

\begin{lstlisting}[language=Python,caption=PyTorch code for neural network architecture used for the SPSD Neural Network surrogate model.,label=listing_two]
import numpy as np
import torch
import torch.nn.functional as F
import torch.nn as nn
torch.set_default_dtype(torch.float64)

class SPD_NN(nn.Module):
    def __init__(self,nHiddenLayers,nNeuronsPerLayer,nFeatures,nOutputs):
        super(SPD_NN, self).__init__()
        forward_list = []
        assert(nFeatures == nOutputs)
        idx = np.tril_indices(nOutputs)
        self.idx = idx
        self.nOutputs = nOutputs
        networkOutputSize = idx[0].size
        self.numHiddenLayers = nHiddenLayers
        self.numLayers = self.numHiddenLayers + 1
        dim = np.zeros(nHiddenLayers+2,dtype='int')
        dim[0] = nFeatures
        for i in range(1,nHiddenLayers+1):
         dim[i] = nNeuronsPerLayer
        dim[-1] = networkOutputSize
        self.dim = dim
        input_dim = dim[0:-1]
        output_dim = dim[1::]

        for i in range(0,nHiddenLayers+1):
          forward_list.append(nn.Linear(input_dim[i], output_dim[i]))

        self.forward_list = nn.ModuleList(forward_list)
        self.activation = F.relu

    def forward(self,x):
      y = x*1.
      for i in range(0,self.numLayers-1):
        y = self.activation(self.forward_list[i](y))
      y = self.forward_list[-1](y)
      K = torch.zeros(y.shape[0],self.nOutputs,self.nOutputs)
      K[:,self.idx[0],self.idx[1]] = y[:]
      KT = torch.transpose(K,2,1)
      K = torch.matmul(K,KT)
      result = torch.einsum('ijk,ik->ij',K,x)
      return result[:,:]
\end{lstlisting}

\section{Proper orthogonal decomposition}\label{appendix:pod}
	Algorithm \ref{alg:pod} provides the algorithm for computing the POD basis
	used in this work. We note that the basis dimension $K$ can be determined
	from the decay of the singular values; for simplicity, we treat it as an
	algorithm input.
\begin{algorithm}
	\caption{Proper orthogonal decomposition (POD)}
\label{alg:pod}
	\textbf{Input}: Snapshot matrix, $\mathbf{S} \in \RR{N \times N_s}$;
	desired basis
	dimension $K$.\\
	\textbf{Output}: POD basis, $\basis \in \RR{N \times K}$, $K \le N_s$. \\
	\textbf{Steps}:
\begin{enumerate}
\item Compute the (thin) singular value decomposition
\begin{equation} 
        \mathbf{S} = \svdLeft
	\svdMid \svdRight,
\end{equation} 
where
$\svdLeft\equiv\left[\svdLeftVecArg{1}\ \cdots\
\svdLeftVecArg{N_s}\right]$.
\item Truncate the left singular vectors such that 
${\basis}_{i} = \svdLeftVecArg{i}$, $i=1,\ldots,K$.
\end{enumerate}
\end{algorithm}

The accuracy of the POD approximation can be bounded by energy contained in the truncated singular values, $\sum_{k = K}^{N_s} \svdMid_{kk}^2$. 

\end{appendix}
\end{document}